\documentclass[]{article}

\usepackage{amssymb,latexsym,amsmath,graphicx,subfigure,color,amsthm,IEEEtrantools}     % Standard packages
\usepackage{stmaryrd}
\usepackage{multirow}

\DeclareMathOperator*{\argmin}{argmin}

\newtheorem{theorem}{Theorem}
\newtheorem{lemma}[theorem]{Lemma}

\begin{document}

\title{Constraint Energy Minimizing Generalized Multiscale Finite Element Method for dual continuum model}
\author{
Siu Wun Cheung\thanks{Department of Mathematics, Texas A\&M University, College Station, TX 77843, USA (\texttt{tonycsw2905@math.tamu.edu})}
\and
Eric T. Chung\thanks{Department of Mathematics, The Chinese University of Hong Kong, Shatin, New Territories, Hong Kong SAR, China (\texttt{tschung@math.cuhk.edu.hk}) }
\and
Yalchin Efendiev\thanks{Department of Mathematics \& Institute for Scientific Computation (ISC), Texas A\&M University,
College Station, Texas, USA (\texttt{efendiev@math.tamu.edu})}
\and
Wing Tat Leung\thanks{
Institute of Computational Engineering and Sciences, University of Texas at Austin, Austin, USA}
\and
Maria Vasilyeva\thanks{Institute for Scientific Computation, Texas A\&M University, College Station, TX 77843 \& Department of Computational
Technologies, North-Eastern Federal University, Yakutsk, Republic of Sakha (Yakutia), Russia, 677980. 
(\texttt{vasilyevadotmdotv@gmail.com})}
}\maketitle

\begin{abstract}
The dual continuum model serves as a powerful tool in the modeling of subsurface applications. It allows a systematic coupling of various components of 
the solutions. The system is of multiscale nature as it involves high heterogeneous and high contrast coefficients. To numerically compute the solutions,
some types of reduced order methods are necessary. We will develop and analyze a novel multiscale method based on
the recent advances in multiscale finite element methods. Our method will compute multiple local multiscale basis functions per coarse region.
The idea is based on some local spectral problems, which are important to identify high contrast channels, and an energy minimization principle.
Using these concepts, we show that the basis functions are localized, even in the presence of high contrast long channels and fractures. 
In addition, we show that the convergence of the method depends only on the coarse mesh size. Finally, we present several
numerical tests to show the performance. 
\end{abstract}

\section{Introduction}

Common in a wide variety of applications related to subsurface formations, one needs to perform 
numerical simulations in domains containing discrete fractures, faults and thin structures. 
The material properties within fractures can have a large difference from the material properties in the background media,
which can also contain highly heterogeneous and high contrast regions. These large contrasts in material properties and 
the complex geometries of the fractures lead to difficulties in traditional numerical simulations due
to the fact that solutions contain various scales and resolving these scales requires huge computational costs.  
Our goal in this paper is to construct and analyze reduced models for such problems. 
In classical upscaling approach, the computational domain is decomposed into coarse blocks, not necessarily resolving scales, 
and effective material property is computed for each coarse block \cite{dur91,weh02}.
To compute effective material properties, some local problems are solved. However, it is known
that one effective coefficient per coarse region is not enough to capture various properties of the solutions,
especially for regions with fractures and high contrast heterogeneities. 
To overcome this drawback, the multi-continuum approaches are used
\cite{arbogast1990derivation,barenblatt1960basic,kazemi1976numerical,pruess1982fluid,warren1963behavior,wu1988multiple},
where several effective medium properties are constructed. 
For example, in flow problems, separate equations for the flow in the background (called matrix)
and the flow within fractures are derived, and these quantities are coupled by some interaction terms. 
The multi-continuum model thus provides a powerful tool for problems for subsurface applications with fractures. 

One important component of our approach is a local fine grid simulation, 
which is typical in many multiscale and numerical upscaling techniques. 
In general, a fine grid simulation involving flow and transport in heterogenous fracture media can be decomposed into two parts (we refer \cite{karimi2016general}
for an overview). First of all, an unstructured fine mesh is needed to model the geometries of the fractures and background heterogeneities. 
Secondly, using the fine mesh, the underlying physical model is discretized. There are in literature a variety of numerical approaches.
For instances, in \cite{baca1984modelling, juanes2002general, karimi2003numerical, kim2000finite}, the standard Galerkin formulation is considered,
in \cite{erhel2009flow, hoteit2008efficient, ma2006mixed, martin2005modeling}, the mixed finite element method is considered,
and in \cite{eikemo2009discontinuous, hoteit2005multicomponent}, the discontinuous Galerkin method is considered. 
Moreover, in \cite{bogdanov2003two, granet2001two, karimi2004efficient, monteagudo2004control, noetinger2015quasi, reichenberger2006mixed}, the finite volume
scheme is investigated. A hybrid scheme combining the finite element method for the pressure equation and the finite volume method for the transport equation
has also been considered \cite{geiger2009black, matthai2007finite, nick2011comparison}.

The reduced model we developed in this paper is motivated by the 
Generalized Multiscale Finite Element Method
(GMsFEM) \cite{egh12,chung2016adaptive,chung2015generalizedwave},
which can be seen as a generalization of the multiscale finite element method (MsFEM). 
We will construct multiscale basis functions that can couple various continua as well as effects of high contrast channels and fractures.
The main idea of GMsFEM is to identify local dominant modes
by the use of local spectral problems defined in some suitable snapshot spaces. 
These ideas are important in identifying influences of high contrast channels and regions,
which are required to be represented individually by separate basis functions. In this regard, 
the GMsFEM shares some similarities with the multi-continuum approaches (see \cite{chung2017coupling}).
The idea of constructing local basis functions using spectral problem has also been used
by various domain decomposition methods \cite{ge09_1reduceddim,kim2015bddc,kim2016bddc}.
We remark that the convergence of the GMsFEM is related to the decay of the eigenvalues
of the local spectral problems \cite{chung2014adaptive}.

%{\bf GMsFEM  convergence in terms of mesh and we use localized ideas}
%It was shown that the GMsFEM's convergence depends on the eigenvalue
%decay \cite{chung2014adaptive}.
It is in general not an easy task to derive a multiscale method with a convergence depends only on the coarse mesh size
and independent of scales and contrast. 
To obtain multiscale methods with mesh dependent convergence, several approaches 
are considered in literature \cite{owhadi2014polyharmonic, maalqvist2014localization, owhadi2017multigrid,hou2017sparse,chung2018constraint,chung2018mixed}.
The theory of GMsFEM motivates the use of local spectral problems to capture the effects of high contrast channels.
This idea is also used in obtaining mesh dependent convergence \cite{hou2017sparse,chung2018constraint,chung2018mixed}.

%{\bf Our method: (1) existing methods close;
%(2) Our method and main difference - contrast..}

In this paper, we will develop and analyze a novel multiscale method
for a dual continuum model with a convergence depends only on the coarse mesh size
and independent of scales and contrast. Our ideas are motivated by the Constraint Energy Minimizing Generalized Multiscale
Finite Element Method (CEM-GMsFEM) \cite{chung2018constraint,chung2018mixed}.
There are two ingredients of our methodology. 
First of all, we will construct a set of local auxiliary multiscale basis functions, as in GMsFEM. 
These functions are dominant eigenfunctions of local spectral problems,
and the number of these functions is the same as the number of high contrast channels and fracture networks.
We emphasize that this is the minimal number of degrees of freedoms required to represent channelized effects. 
We also remark that these eigenfunctions are crucial in the construction of localized basis functions. 
The second key component is multiscale basis functions. These functions are obtained
by minimizing an energy functional subject to certain constraints. These constraints
are formulated using the auxiliary functions with the purpose of obtaining localized multiscale basis functions. 
In particular, for each of the auxiliary function, the constraints require the minimizer of the energy functional
is orthogonal, in a weighted $L^2$ sense, to all other auxiliary functions except the selected one. For the selected auxiliary functions,
the  constraints require the minimizer of the energy functional to satisfy a normalized condition. 
Combining the effects of auxiliary functions and energy minimization, we show that the minimizer of the energy functional
has exponential decay property, and is very small outside an oversampling region obtained by the support of the selected auxiliary function. 
Moreover, the resulting multiscale method obtained by a Galerkin formulation has a mesh dependent convergence rate.
We remark that one can also perform adaptivity as in \cite{chung2014adaptive,chung2015residual,chung2018fast}.

The paper is organized as follows. In Section~\ref{sec:model}, we will introduce the dual continuum model. 
Our multiscale method will be presented in Section~\ref{sec:method} and analyzed in Section~\ref{sec:analysis}.
In Section~\ref{sec:numerical}, we will present some numerical tests. The paper ends with a conclusion in Section~\ref{sec:conclusions}.

\section{Dual continuum Model}\label{sec:model}
We consider the following dual continuum model
\begin{equation}
\begin{split}
c_1 \dfrac{\partial p_1}{\partial t} - \text{div}(\kappa_1 \nabla p_1) + \rho \sigma(p_1 - p_2) = \rho f_1, \\
c_2 \dfrac{\partial p_2}{\partial t} - \text{div}(\kappa_2 \nabla p_2) - \rho \sigma(p_1 - p_2) = \rho f_2,
\end{split}
\label{eq:dc}
\end{equation}
in a computational domain $\Omega \subset \mathbb{R}^d$. The domain $\Omega$ is divided into the fracture and the matrix region
\begin{equation}
\Omega = D_m \oplus_i d_i D_{f,i},
\end{equation}
where $m$ and $f$ represent the matrix and the fracture regions. $d_i$ denotes the aperture of the $i$-th fracture and $i$ is the index of the fractures. We denote by $\kappa_i$ the permeability of the $i$-th fracture. 
The continua are coupled through a mass exchange in the last term on the left hand side of \eqref{eq:dc}. 
$D_m$ is a two-dimensional domain and $D_{f,i}$ is a one-dimensional domain. We prescribe the initial condition $p_i(0,\cdot) = p_i^0$ in $\Omega$ and the boundary condition $p_i(t, \cdot) = 0$ on $\partial \Omega$ for $t > 0$.
Here, we assume the permeability fields are uniformly bounded, i.e. 
\begin{equation}
0 < \underline{\kappa} \leq \kappa_i(x), \kappa_{l,i}(x) \leq \overline{\kappa} \text{ for } x \in \Omega.
\end{equation}

Let $V = [H^1_0(\Omega)]^2$. 
Also, for a subdomain $D \subset \Omega$, we denote the restriction of $V$ on $D$ by $V(D)$, and
the subspace of $V(D)$ with zero trace on $\partial D$ by $V_0(D)$.
The weak formulation of \eqref{eq:dc} then reads:
find $p = (p_1, p_2)$ such that $p(t, \cdot) \in V$ and
\begin{equation}
\begin{split}
c \left(\dfrac{\partial p}{\partial t}, v \right) + a_Q(p, v) = (f,v),
\end{split}
\label{eq:sol_weak}
\end{equation}
for all $v = (v_1, v_2)$ with $v(t, \cdot) \in V$.
The bilinear forms are defined as:
\begin{equation}
\begin{split}
c_i(u_i, v_i) & = \int_{D_m} c_i u_i v_i \, dx + \sum_l \int_{D_{f,l}} c_{l,i} u_i v_i \, ds,\\
c(u,v) & = \sum_i c_i(u_i,v_i), \\
a_i(p_i, v_i) & = \int_{D_m} \kappa_i \nabla p_i \cdot \nabla v_i \, dx + \sum_l \int_{D_{f,l}} \kappa_{l,i} \nabla_f p_i \cdot \nabla_f v_i \, ds, \\
a(p,v) & = \sum_i a_i(p_i, v_i), \\
q(p,v) & = \sum_i \sum_l \rho \sigma \int_\Omega (p_i - p_l) v_i \, dx,\\
a_Q(p,v) & = a(p,v) + q(p,v).
%(f,v) & =  \int_{D_m} f_i v_i \, dx + \sum_l \int_{D_{f,l}} f_i v_i \, ds.
\end{split}
\label{eq:bilinear}
\end{equation}

\section{Method description}\label{sec:method}

\begin{figure}[ht!]
\centering
\includegraphics[width=0.5\linewidth]{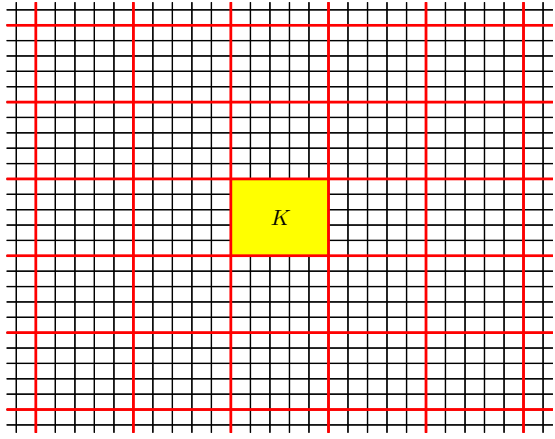}
\caption{An illustration of the fine grid and the coarse grid and a coarse element.}
\label{fig:mesh}
\end{figure}

In this section, we will describe the details of our proposed method. 
To start with, we introduce the notions of coarse and fine meshes. 
We start with a usual partition $\mathcal{T}^H$ of $\Omega$ into finite elements, 
which does not necessarily resolve any multiscale features. 
The partition $\mathcal{T}^H$ is called a coarse grid and 
a generic element $K$ in the partition $\mathcal{T}^H$ is called a coarse element. 
Moreover, $H > 0$ is called the coarse mesh size.
We let $N_c$ be the number of coarse grid nodes and 
$N$ be the number of coarse elements. 
We also denote the collection of all coarse grid edges by $\mathcal{E}^H$.
We perform a refinement of $\mathcal{T}^H$ to obtain a fine grid $\mathcal{T}^h$, 
where $h > 0$ is called the fine mesh size. 
It is assumed that the fine grid is sufficiently fine to resolve the solution. 
An illustration of the fine grid and the coarse grid and a coarse element are shown in Figure~\ref{fig:mesh}.

We define local bilinear forms on a coarse element $K_j$ by:
\begin{equation}
\begin{split}
a^{(j)}_i(p_i, v_i) & = \int_{K_j} \kappa_i \nabla p_i \cdot \nabla v_i \, dx + \sum_l \int_{D_{f,l} \cap K_j} \kappa_{l,i} \nabla_f p_i \cdot \nabla_f v_i \, ds, \\
a^{(j)}(p,v) & = \sum_i a^{(j)}_i(p_i, v_i), \\
q^{(j)}(p,v) & = \sum_i \sum_l \rho \sigma \int_{K_j} (p_i - p_l) v_i \, dx,\\
a^{(j)}_Q(p,v) & = a^{(j)}(p,v) + q^{(j)}(p,v), \\
s^{(j)}_i(p_i, v_i) & = \int_{K_j} \widetilde{\kappa}_i p_i v_i \, dx + \sum_l \int_{D_{f,l} \cap K_j} \widetilde{\kappa}_{l,i} p_i v_i \, ds, \\
s^{(j)}(p,v) & = \sum_i s^{(j)}_i(p_i, v_i),
\end{split}
\label{eq:bilinear_loc}
\end{equation}
where $\widetilde{\kappa}_i = \kappa_i \sum_{k=1}^{N_c} \vert \nabla \chi_k \vert^2$, $\widetilde{\kappa}_{l,i} = \kappa_{l,i} \sum_{k=1}^{N_c} \vert \nabla_f \chi_k \vert^2$, and $\{\chi_k\}$ is a set of bilinear partition of unity functions for the coarse grid partition of the domain $\Omega$. We also define the bilinear form $s$ by:
\begin{equation}
s(p,v) = \sum_j s^{(j)}(p,v).
\end{equation}

Next, we will use the concept of GMsFEM to construct our auxiliary multiscale basis functions.
The auxiliary basis functions are coupled, and defined by a spectral problem, 
which is to find a real number $\lambda_k^{(j)}$ 
and a function $\phi_k^{(j)} \in V(K_j)$ such that
\begin{equation}
a_Q^{(j)}(\phi_k^{(j)}, v) = \lambda_k^{(j)} s^{(j)}(\phi_k^{(j)}, v) \text{ for all } v \in V(K_j).
\label{eq:spectral_prob}
\end{equation}
We let $\lambda_k^{(j)}$ be the eigenvalues of \eqref{eq:spectral_prob} 
arranged in ascending order in $k$, and use the first $L_j$ eigenfunctions 
to construct our local auxiliary multiscale space
\begin{equation}
V_{aux}^{(j)} = \text{span} \{ \phi_k^{(j)}: 1 \leq k \leq L_j\}.
\end{equation}
The global auxiliary multiscale space $V_{aux}$ is then defined as 
the sum of these local auxiliary multiscale spaces
\begin{equation}
V_{aux} = \oplus_{j=1}^N V_{aux}^{(j)}.
\end{equation}

Before we move on to discuss the construction of multiscale basis functions, 
we introduce some tools which will be used to describe our method and analyze the convergence.
We first introduce the notion of $\phi$-orthogonality. 
In a coarse block $K_j$, given an auxiliary basis function $\phi_k^{(j)} \in V_{aux}$,
we say that $\psi \in V$ is $\phi_k^{(j)}$-orthogonal if
\begin{equation}
s \left(\psi, \phi_{k'}^{(j')}\right) = \delta_{j,j'} \delta_{k,k'}.
\end{equation}
We also introduce a projection operator $\pi: [L^2(\Omega)]^2 \to V_{aux}$ by 
$\pi = \sum_{j=1}^N \pi_j$, where $\pi_j: [L^2(K_j)]^2 \to V_{aux}$ is given by
\begin{equation}
\pi_j(v) = \sum_{k=1}^{L_j} \dfrac{s^{(j)}(v,\phi_k^{(j)})}{s^{(j)}(\phi_k^{(j)},\phi_k^{(j)})} \phi_k^{(j)} \text{ for all } v \in [L^2(K_j)]^2.
\end{equation}

Next, we construct our global multiscale basis functions.
The global multiscale basis function $\psi_{j}^{(i)} \in V$ is defined as the solution of 
the following constrained energy minimization problem
\begin{equation}
\psi_{k}^{(j)} = \argmin \left\{ a_Q(\psi, \psi) : 
\psi \in V \text{ is } \phi_k^{(j)} \text{-orthogonal}\right\}.
\label{eq:min1_glo}
\end{equation}
The minimization problem \eqref{eq:min1_glo} is equivalent to the following variational problem: 
find $\psi_{k}^{(j)} \in V$ and $\mu_{k}^{(j)} \in V_{aux}^{(j)}$ such that
\begin{equation}
\begin{split}
a_Q(\psi_{k}^{(j)}, w) + s^{(j)}(w, \mu_{k}^{(j)}) & = 0 \text{ for all } w \in V, \\
s^{(j)}(\psi_{k}^{(j)} - \phi_k^{(j)}, \nu) & = 0 \text{ for all } \nu \in V_{aux}^{(j)}.
\end{split}
\label{eq:var1_glo}
\end{equation}

\begin{figure}[ht!]
\centering
\includegraphics[width=0.5\linewidth]{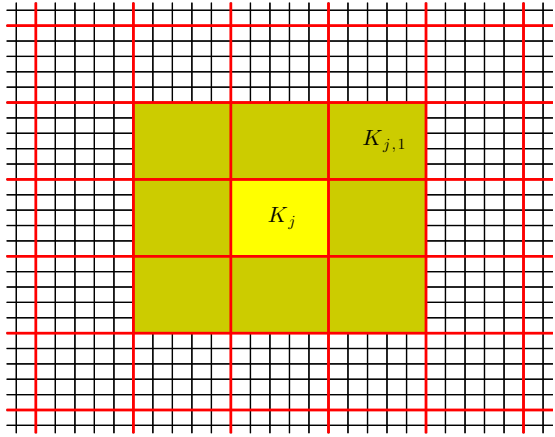}
\caption{An illustration of an oversampled domain formed by enlarging $K_j$ with $1$ coarse grid layer.}
\label{fig:oversample}
\end{figure}

Motivated by the construction of global multiscale basis functions,
we define our localized multiscale basis functions.
For each element $K_j$, an oversampled domain formed by 
enlarging the coarse grid block $K_i$ by $m$ coarse grid layers. 
An illustration of an oversampled domain is shown in Figure~\ref{fig:oversample}.
The localized multiscale basis function $\psi_{k,{ms}}^{(j)} \in V_0(K_{j,m})$ is defined as the solution of 
the following constrained energy minimization problem
\begin{equation}
\psi_{k,{ms}}^{(j)} = \argmin \left\{ a_Q(\psi, \psi) : 
\psi \in V_0(K_{j,m}) \text{ is } \phi_k^{(j)} \text{-orthogonal}\right\}.
\label{eq:min1}
\end{equation}
The minimization problem \eqref{eq:min1} is equivalent to the following variational problem: 
find $\psi_{k,{ms}}^{(j)} \in V_0(K_{j,m})$ and $\mu_{k,ms}^{(j)} \in V_{aux}^{(j)}$ such that
\begin{equation}
\begin{split}
a_Q(\psi_{k,{ms}}^{(j)}, w) + s^{(j)}(w, \mu_{k,ms}^{(j)}) & = 0 \text{ for all } w \in V_0(K_{j,m}), \\
s^{(j)}(\psi_{k,{ms}}^{(j)} - \phi_k^{(j)}, \nu) & = 0 \text{ for all } \nu \in V_{aux}^{(j)}.
\end{split}
\label{eq:var1}
\end{equation}

We use the localized multiscale basis functions to construct
the multiscale finite element space, which is defined as
\begin{equation}
V_{ms} = \text{span}  \{\psi_{k,{ms}}^{(j)} : 1 \leq k \leq L_j, 1 \leq j \leq N \}.
\end{equation}
The multiscale solution is then given by: find $p_{ms} = (p_{ms,1}, p_{ms,2})$ with $p_{ms}(t,\cdot) \in V_{ms}$ such that for all $v = (v_1, v_2)$ with $v(t, \cdot) \in V_{ms}$,
\begin{equation}
c \left(\dfrac{\partial p_{ms}}{\partial t}, v \right) + a_Q(p_{ms}, v) = (f,v).
\label{eq:sol_ms}
\end{equation}

\section{Convergence Analysis}\label{sec:analysis}

In this section, we will analyze the proposed method.
First, we define the following norms and semi-norms on $V$:
\begin{equation}
\begin{split}
\| p \|_c^2 & = c(p,p), \\ 
\| p \|_a^2 & = a(p,p), \\
\vert p \vert_q^2 & = q(p,p), \\
\| p \|_{a_Q}^2 & = a_Q(p,p), \\
\| p \|_s^2 & = s(p,p).
\end{split}
\end{equation}
For a subdomain $D = \bigcup_{j \in J} K_j$ composed by a union of coarse grid blocks, we also define the following local norms and semi-norms on $V$:
\begin{equation}
\begin{split}
\| p \|_{a(D)}^2 & = \sum_{j \in J} a^{(j)}(p,p), \\
\vert p \vert_{q(D)}^2 & = \sum_{j \in J} q^{(j)}(p,p), \\
\| p \|_{a_Q(D)}^2  & = \sum_{j \in J} a_Q^{(j)}(p,p), \\
\| p \|_{s(D)}^2 & = \sum_{j \in J} s^{(j)}(p,p).
\end{split}
\end{equation}
The flow of our analysis goes as follows.
First, we prove the convergence using the global multiscale basis functions.
With the global multiscale basis functions constructed, 
the global multiscale finite element space is defined by
\begin{equation}
V_{glo} = \text{span}  \{\psi_{k}^{(j)} : 1 \leq k \leq L_j, 1 \leq j \leq N \},
\end{equation}
and an approximated solution $p_{glo} = (p_{glo,1}, p_{glo,2})$, 
where $p_{glo}(t,\cdot) \in V_{glo}$, is given by
\begin{equation}
c \left(\dfrac{\partial p_{glo}}{\partial t}, v \right) + a_Q(p_{glo}, v) = (f,v),
\label{eq:sol_glo}
\end{equation}
for all $v = (v_1, v_2)$ with $v(t,\cdot) \in V_{glo}$.
Next, we give an estimate of the difference between 
the global multiscale functions $\psi_k^{(j)}$ and 
the local multiscale basis functions $\psi_{k,ms}^{(j)}$, 
in order to show that using the multiscale solution $p_{ms}$ 
provides similar convergence results as the global solution $p_{glo}$.
For this purpose, we denote the kernel of 
the projection operator $\pi$ by $\widetilde{V}$. 
Then, for any $\psi_k^{(j)} \in V_{glo}$, we have 
\begin{equation}
a_Q(\psi_k^{(j)}, w) = 0 \text{ for all } w \in \widetilde{V},
\end{equation}
which implies $\widetilde{V} \subseteq V_{glo}^\perp$, 
where $V_{glo}^\perp$ is the orthogonal complement of 
$V_{glo}$ with respect to the inner product $a_Q$. 
Moreover, since $\text{dim}(V_{glo}) = \text{dim}(V_{aux})$, 
we have $\widetilde{V} = V_{glo}^\perp$ and 
$V = V_{glo} \oplus \widetilde{V}$.

In addition, we introduce some operators which will be used in our analysis, 
namely $R_{glo}: V \to V_{glo}$ given by: 
for any $u \in V$, the image $R_{glo}u \in V_{glo}$ is defined by
\begin{equation}
a_Q(R_{glo}u , v) = a_Q(u,v) \text{ for all } v \in V_{glo},
\label{eq:elliptic_proj}
\end{equation}
and similarly, $R_{ms}: V \to V_{ms}$ given by: 
for any $u \in V$, the image $R_{ms}u \in V_{ms}$ is defined by
\begin{equation}
a_Q(R_{ms}u , v) = a_Q(u,v) \text{ for all } v \in V_{ms}.
\label{eq:elliptic_proj_ms}
\end{equation}
We also define $\mathcal{C}: V \to V$ given by:
for any $u \in V$, the image $\mathcal{C}u \in V$ is defined by
\begin{equation}
(\mathcal{C}u , v) = c(u,v) \text{ for all } v \in V.
\end{equation}
Moreover, the operator $\mathcal{A}: D(\mathcal{A}) \to [L^2(\Omega)]^2$ is defined on a subspace $D(\mathcal{A}) \subset V$ by:
for any $u \in D(\mathcal{A})$, the image $\mathcal{A}u \in [L^2(\Omega)]^2$ is defined by
\begin{equation}
(\mathcal{A}u , v) = a_Q(u,v) \text{ for all } v \in V.
\end{equation}

We will first show the projection operator $R_{glo}$ onto global multiscale finite element space 
has a good approximation property with respect to the $a_Q$-norm and $L^2$-norm.
\begin{lemma}
Let $u \in D(\mathcal{A})$. Then we have $u - R_{glo} u \in \widetilde{V}$ and
\begin{equation}
\| u - R_{glo} u \|_{a_Q} \leq CH \underline{\kappa}^{-\frac{1}{2}} \Lambda^{-\frac{1}{2}} \| \mathcal{A} u \|_{[L^2(\Omega)]^2},
\label{eq:a_approx_glo}
\end{equation}
and
\begin{equation}
\| u - R_{glo} u \|_{[L^2(\Omega)]^2} \leq CH^2 \underline{\kappa}^{-1} \Lambda^{-1} \| \mathcal{A} u \|_{[L^2(\Omega)]^2},
%\| u - R_{glo} u \|_c \leq CH^2 \Lambda^{-1} \| \mathcal{A} u \|_{[L^2(\Omega)]^2},
\label{eq:c_approx_glo}
\end{equation}
where
\begin{equation}
\Lambda = \min_{1 \leq j \leq N} \lambda_{L_j+1}^{(j)}.
\end{equation}
\begin{proof}
From \eqref{eq:elliptic_proj}, we see that
$u - R_{glo} u \in V_{glo}^\perp = \widetilde{V}$. 
Taking $v = R_{glo} u \in V_{glo}$ in \eqref{eq:elliptic_proj}, we have
\begin{equation}
\begin{split}
a_Q(u - R_{glo}u, R_{glo} u) = 0.
\end{split}
\end{equation}
Therefore, we have
\begin{equation}
\begin{split}
\|u - R_{glo} u\|_{a_Q}^2 
& = a_Q(u - R_{glo}u, u - R_{glo} u) \\
& = a_Q(u - R_{glo}u, u) \\
& = a_Q(u, u - R_{glo} u) \\
& = (\mathcal{A}u, u - R_{glo} u) \\
& \leq \| \widetilde{\kappa}^{-\frac{1}{2}} \mathcal{A} u\|_{[L^2(\Omega)]^2} \| u - R_{glo} u \|_s,
\end{split}
\end{equation}
where $\widetilde{\kappa}(x) = \min \{ \widetilde{\kappa}_i(x), \widetilde{\kappa}_{l,i}(x)\}$.
Since $u - R_{glo} u \in \widetilde{V}$, we have $\pi_j(u - R_{glo}u) = 0$ for all $j = 1,2,\ldots,N$ and
\begin{equation}
\begin{split}
\| u - R_{glo} u \|_s^2
& = \sum_{j=1}^N \| u - R_{glo} u \|_{s(K_j)}^2 \\
& = \sum_{j=1}^N \| (I - \pi_j) (u - R_{glo} u) \|_{s(K_j)}^2.
\end{split}
\end{equation}
By the orthogonality of the eigenfunctions $\phi_k^{(j)}$, we have
\begin{equation}
\begin{split}
\sum_{j=1}^N \| (I - \pi_j) (u - R_{glo} u) \|_{s(K_j)}^2 
\leq \dfrac{1}{\Lambda} \sum_{j=1}^N \| u - R_{glo} u \|_{a_Q(K_j)}^2
\leq \dfrac{1}{\Lambda} \| u - R_{glo} u \|_{a_Q}^2.
\end{split}
\end{equation}
Finally, using the fact that $\vert \nabla \chi_k \vert = O(H^{-1})$, we obtain the first estimate \eqref{eq:a_approx_glo}. 

For the second estimate \eqref{eq:c_approx_glo}, we use a duality argument. Define $w \in V$ by
\begin{equation}
a_Q(w, v) = (u - R_{glo} u, v) \text{ for all } v \in V.
\end{equation}
Then we have
\begin{equation}
\begin{split}
\| u - R_{glo} u \|_{[L^2(\Omega)]^2}^2 = (u - R_{glo} u, u - R_{glo} u) = a_Q(w, u - R_{glo} u).
\end{split}
\end{equation}
Taking $v = R_{glo} w \in V_{glo}$ in \eqref{eq:elliptic_proj}, we have
\begin{equation}
\begin{split}
a_Q(u - R_{glo}u, R_{glo} w) = 0.
\end{split}
\end{equation}
Note that $w \in D(\mathcal{A})$ and $\mathcal{A}w = u - R_{glo} u$. Hence
\begin{equation}
\begin{split}
\| u - R_{glo} u \|_{[L^2(\Omega)]^2}^2
& = a_Q(w - R_{glo} w, u - R_{glo} u) \\
& \leq \| w - R_{glo} w \|_{a_Q} \| u - R_{glo} u \|_{a_Q} \\
& \leq \left(CH \underline{\kappa}^{-\frac{1}{2}} \Lambda^{-\frac{1}{2}} \| \mathcal{A} w \|_{[L^2(\Omega)]^2}\right) \left( CH \underline{\kappa}^{-\frac{1}{2}}\Lambda^{-\frac{1}{2}} \| \mathcal{A} u \|_{[L^2(\Omega)]^2}\right) \\
& \leq CH^2 \underline{\kappa}^{-1}\Lambda^{-1} \| u - R_{glo} u \|_{[L^2(\Omega)]^2} \| \mathcal{A} u \|_{[L^2(\Omega)]^2}.
\end{split}
\end{equation}
\end{proof}
\label{thm:elliptic-glo}
\end{lemma}

We are now going to prove the global basis functions are localizable. 
For each coarse block $K$, we define $B$ to be a bubble function with $B(x) > 0$ for all $x \in \text{int}(K)$ 
and $B(x) = 0$ for all $x \in \partial K$. 
We will take $B = \prod_j \chi_j^{ms}$ where the product is taken over all vertices 
$j$ on the boundary of $K$. Using the bubble function, we define the constant
\begin{equation}
C_\pi = \sup_{K \in \mathcal{T}^H, \nu \in V_{aux}} \dfrac{s(\nu,\nu)}{s(B\nu, \nu)}.
\end{equation}
We also define
\begin{equation}
\lambda_{max} = \max_{1 \leq j \leq N} \max_{1 \leq k \leq L_j} \lambda_k^{(j)}.
\end{equation}

\begin{lemma}
\label{lemma2}
For all $v_{aux} \in V_{aux}$, there exists a function $v \in V$ such that
\begin{equation}
\pi(v) = v_{aux}, \quad \| v \|_{a_Q}^2 \leq D \| v_{aux} \|_s^2, \quad \text{supp}(v) \subset \text{supp}(v_{aux}).
\end{equation}
We write $D = 2(1+2C_p^2 \rho \sigma \underline{\kappa}^{-1})(C_\mathcal{T} + \lambda_{max}^2)$, where $C_\mathcal{T}$ is the square of the maximum number of vertices over all coarse elements, and $C_p$ is a Poincar\'{e} constant.
\begin{proof}
Let $v_{aux} \in V_{aux}^{(j)}$ with $\|v_{aux}\|_{s(K_j)} = 1$. We consider the following minimization problem defined on a coarse block $K_j$.
\begin{equation}
v = \text{argmin}\left\{ a_Q(\psi, \psi) \: : \: \psi \in V_0(K_j), \quad s^{(j)}(\psi, \nu) = s^{(j)}(v_{aux}, \nu) \text{ for all } \nu \in V_{aux}^{(j)} \right\}.
\label{eq:min_lemma2}
\end{equation}
We will show that the minimization problem \eqref{eq:min_lemma2} has a unique solution.
First, we note that the minimization problem \eqref{eq:min_lemma2} is equivalent to the following variational problem: find $v \in V_0(K_j)$ and $\mu \in V_{aux}^{(j)}$ such that
\begin{equation}
\begin{split}
a_Q^{(j)}(v,w) + s^{(j)}(w,\mu) & = 0 \text{ for all } w \in V_0(K_j), \\
s^{(j)}(v - v_{aux},\nu) & = 0 \text{ for all } \nu \in V_{aux}^{(j)}.
\end{split}
\label{eq:var_lemma2}
\end{equation}
The well-posedness of \eqref{eq:var_lemma2} is equivalent to the existence of $v \in V_0(K_j)$ such that
\begin{equation}
s^{(j)}(v, v_{aux}) \geq C \| v_{aux} \|_{s(K_j)}^2, \quad \| v \|_{a_Q(K_j)} \leq C \|v_{aux}\|_{s(K_j)}, 
\end{equation}
where $C$ is a constant to be determined. Now, we take $v = Bv_{aux} \in V_0(K_j)$. Then we have
\begin{equation}
s^{(j)}(v, v_{aux}) = s^{(j)}(Bv_{aux}, v_{aux}) \geq C_\pi^{-1} s \| v_{aux} \|_{s(K_j)}^2.
\end{equation}
On the other hand, since $\nabla v_i = \nabla (B v_{aux,i}) = v_{aux,i} \nabla B + B \nabla v_{aux,i}$, $\vert B \vert \leq 1$ and $\vert \nabla B \vert^2 \leq C_{\mathcal{T}} \sum_k \vert \nabla \chi_k^{ms} \vert^2$, we have
\begin{equation}
\| v \|_{a(K_j)}^2 \leq 2(C_\mathcal{T} \| v_{aux} \|_{s(K_j)}^2 + \| v_{aux} \|_{a_Q(K_j)}^2).
\end{equation}
By the spectral problem \eqref{eq:spectral_prob}, we have
\begin{equation}
\| v_{aux} \|_{a_Q(K_j)} \leq \max_{1 \leq k \leq L_j} \lambda_k^{(j)}\| v_{aux} \|_{s(K_j)}.
\end{equation}
Moreover, by Poincar\'{e} inequality, we have
\begin{equation}
\vert v \vert_q^2 \leq 2 \rho \sigma \| v \|_{L^2(K_j)}^2 \leq 2C_p^2 \rho \sigma \underline{\kappa}^{-1} \| v \|_{a(K_j)}^2.
\end{equation}
Combining these estimates, we have
\begin{equation}
\| v \|_{a_Q(K_j)}^2 \leq (1+2C_p^2 \rho \sigma \underline{\kappa}^{-1}) \| v \|_{a(K_j)}^2 
\leq 2(1+2C_p^2 \rho \sigma \underline{\kappa}^{-1})(C_\mathcal{T} + \lambda_{max}^2) \| v_{aux} \|_{s(K_j)}^2.
\end{equation}
This shows that the minimization problem \eqref{eq:min_lemma2} has a unique solution $v \in V_0(K_j)$, which satisfies our desired properties. 
\end{proof}
\end{lemma}

Here, we make a remark that we can assume $D \geq 1$ without loss of generality.

In order to estimate the difference between the global basis functions and localized basis functions, 
we need the notion of a cutoff function with respect to the oversampling regions. 
For each coarse grid $K_j$ and $M > m$, we define $\chi_j^{M,m} \in \text{span}\{ \chi_k^{ms} \}$ 
such that $0 \leq \chi_j^{M,m} \leq 1$ and $\chi_j^{M,m} = 1$ on the inner region $K_{j,m}$ and 
$\chi_j^{M,m} = 0$ outside the region $K_{j,M}$.

The following lemma shows that our multiscale basis functions have a decay property. In particular, the
global basis functions are small outside an oversampled region specified in the lemma, 
which is important in localizing the multiscale basis functions.

\begin{lemma}
\label{lemma3}
Given $\phi_k^{(j)} \in V_{aux}^{(j)}$ and an oversampling region $K_{j,m}$ with number of layers $m \geq 2$. Let $\psi_{k,ms}^{(j)}$ be a localized multiscale basis function defined on $K_{j,m}$ given by \eqref{eq:min1}, and $\psi_k^{(j)}$ be the corresponding global basis function given by \eqref{eq:min1_glo}. Then we have
\begin{equation}
\| \psi_k^{(j)} - \psi_{k,ms}^{(j)} \|_{a_Q}^2 \leq E \| \phi_k^{(j)} \|_{s(K_j)}^2,
\end{equation}
where $E = 24D^2(1+\Lambda^{-1})\left(1+\dfrac{\Lambda^\frac{1}{2}}{3D^\frac{1}{2}}\right)^{1-m}$.
\begin{proof}
By Lemma~\ref{lemma2}, there exists $\widetilde{\phi}_k^{(j)} \in V$ such that
\begin{equation}
\pi(\widetilde{\phi}_k^{(j)}) = \phi_k^{(j)}, \quad \| \widetilde{\phi}_k^{(j)} \|_{a_Q}^2 \leq D \| \phi_k^{(j)} \|_s^2, \quad \text{supp}(\widetilde{\phi}_k^{(j)}) \subset K_j.
\label{eq:lemma3.0}
\end{equation}
We take $\eta = \psi_k^{(j)} - \widetilde{\phi}_k^{(j)} \in V$
and $\zeta = \widetilde{\phi}_k^{(j)} - \psi_{k,ms}^{(j)} \in V_0(K_{j,m})$. 
Then $\pi(\eta) = \pi(\zeta) = 0$ and hence $\eta, \zeta \in \widetilde{V}$. 
Again, by Lemma~\ref{lemma2}, there exists $\beta \in V$ such that
\begin{equation}
\pi(\beta) = \pi(\chi_j^{m,m-1} \eta), \quad 
\| \beta \|_{a_Q}^2 \leq D \| \pi(\chi_j^{m,m-1} \eta) \|_s^2, \quad 
\text{supp}(\beta) \subset K_{j,m} \setminus K_{j,m-1}.
\label{eq:lemma3.1}
\end{equation}
Take $\tau = \beta - \chi_j^{m,m-1} \eta \in V_0(K_{j,m})$. 
Again, $\pi(\tau) = 0$ and hence $\tau \in \widetilde{V}$.
Now, by the variational problems \eqref{eq:var1_glo} and \eqref{eq:var1}, we have
\begin{equation}
a_Q(\psi_k^{(j)} - \psi_{k,ms}^{(j)}, w) + s^{(j)}(w, \mu_k^{(j)} - \mu_{k,ms}^{(j)}) = 0 \text{ for all } w \in V_0(K_{j,m}).
\end{equation}
Taking $w = \tau - \zeta \in V_0(K_{j,m})$ and using the fact that $\tau - \zeta \in \widetilde{V}$, we have
\begin{equation}
a_Q(\psi_k^{(j)} - \psi_{k,ms}^{(j)}, \tau - \zeta) = 0,
\end{equation}
which implies
\begin{equation}
\begin{split}
\|\psi_k^{(j)} - \psi_{k,ms}^{(j)}\|_{a_Q}^2 
& = a_Q(\psi_k^{(j)} - \psi_{k,ms}^{(j)}, \psi_k^{(j)} - \psi_{k,ms}^{(j)}) \\
& = a_Q(\psi_k^{(j)} - \psi_{k,ms}^{(j)}, \eta+\zeta) \\
& = a_Q(\psi_k^{(j)} - \psi_{k,ms}^{(j)}, \eta+\tau) \\
& \leq \| \psi_k^{(j)} - \psi_{k,ms}^{(j)} \|_{a_Q} \| \eta+\tau \|_{a_Q}.
\end{split}
\end{equation}
Therefore, we have
\begin{equation}
\begin{split}
\|\psi_k^{(j)} - \psi_{k,ms}^{(j)}\|_{a_Q}^2 
& \leq \| \eta+\tau \|_{a_Q}^2 \\
& = \| (1 - \chi_j^{m,m-1}) \eta + \beta \|_{a_Q}^2 \\
& \leq 2\left(\| (1 - \chi_j^{m,m-1}) \eta\|_{a_Q}^2 + \| \beta \|_{a_Q}^2\right).
\end{split}
\label{eq:lemma3.2}
\end{equation}
For the first term on the right hand side of \eqref{eq:lemma3.2}, since 
$\nabla \left((1 - \chi_j^{m,m-1}) \eta_i\right) = (1 - \chi_j^{m,m-1}) \nabla \eta_i - \eta_i \nabla \chi_j^{m,m-1}$ 
and $\vert 1 - \chi_j^{m,m-1} \vert \leq 1$, we have
\begin{equation}
\| (1 - \chi_j^{m,m-1}) \eta_i\|_{a_i}^2
\leq 2 \left( \| \eta_i \|^2_{a_i(\Omega \setminus K_{j,m-1})} + \| \eta_i \|^2_{s_i(\Omega \setminus K_{j,m-1})} \right).
\end{equation}
On the other hand, we have
\begin{equation}
\vert (1 - \chi_j^{m,m-1}) \eta \vert_q^2 \leq \vert \eta \vert_{q(\Omega \setminus  K_{j,m-1})}^2.
\end{equation}
For the second term on the right hand side of \eqref{eq:lemma3.2}, using \eqref{eq:lemma3.1} and $\vert \chi_j^{m,m-1} \vert \leq 1$, we have
\begin{equation}
\begin{split}
 \| \beta \|_{a_Q}^2
\leq D \| \pi(\chi_j^{m,m-1} \eta) \|_s^2 
\leq D \| \chi_j^{m,m-1} \eta \|_s^2 
\leq D \| \eta \|_{s(\Omega \setminus K_{j,m-1})}^2.
\end{split}
\end{equation}
Since $\eta \in \widetilde{V}$, by the spectral problem \eqref{eq:spectral_prob}, we obtain
\begin{equation}
\| \eta \|^2_{s(\Omega \setminus K_{j,m-1})} \leq \Lambda^{-1} \| \eta \|^2_{a_Q(\Omega \setminus K_{j,m-1})}.
\end{equation}
Combining these estimates, we have
\begin{equation}
\|\psi_k^{(j)} - \psi_{k,ms}^{(j)}\|_{a_Q}^2 
\leq (4+4\Lambda^{-1}+2D\Lambda^{-1}) \| \eta \|^2_{a_Q(\Omega \setminus K_{j,m-1})}
\leq 6D(1+\Lambda^{-1}) \| \eta \|^2_{a_Q(\Omega \setminus K_{j,m-1})}.
\label{eq:lemma3.3}
\end{equation}
Next, we will prove a recursive estimate for $\| \eta \|^2_{a_Q(\Omega \setminus K_{j,m-1})}$. 
We take $\xi = 1 - \chi_j^{m-1,m-2}$. 
Then $\xi = 1$ in $\Omega \setminus K_{j,m-1}$ and $0 \leq \xi \leq 1$. Hence, using $\nabla(\xi^2 \eta_i) = \xi^2 \nabla \eta_i + 2 \xi \eta_i\nabla \xi$, we have
\begin{equation}
\| \eta \|_{a_Q(\Omega \setminus K_{j,m-1})}^2 
\leq \| \xi \eta \|_{a_Q}^2 
\leq a_Q(\eta, \xi^2 \eta)
+ 2 \| \xi\eta \|_s \| \eta \|_{a_Q(K_{j,m-1} \setminus K_{j,m-2})}.
\label{eq:lemma3.4}
\end{equation}
We will estimate the first term on the right hand side of \eqref{eq:lemma3.4}. 
First, we note that, for any coarse element $K_{j'} \subset \Omega \setminus K_{j,m-1}$, 
since $\xi = 1$ in $K_{j'}$ and $\eta \in \widetilde{V}$, we have
\begin{equation}
s\left(\xi^2 \eta, \phi_{k'}^{(j')}\right) = s\left(\eta, \phi_{k'}^{(j')}\right) = 0 \text{ for all } k' = 1,2,\ldots,L_{j'}.
\end{equation}
On the other hand, for any coarse element $K_{j'} \subset K_{j,m-2}$, since $\xi = 0$ in $K_{j, m-2}$, we have
\begin{equation}
s\left(\xi^2 \eta, \phi_{k'}^{(j')}\right) = 0 \text{ for all } k' = 1,2,\ldots,L_{j'}.
\end{equation}
Therefore, $\text{supp}(\xi^2 \eta) \subset  K_{j,m-1} \setminus K_{j,m-2}$.
By Lemma~\ref{lemma2}, there exists $\gamma \in V$ such that
\begin{equation}
\pi(\gamma) = \pi(\xi^2 \eta), \quad 
\| \gamma \|_{a_Q}^2 \leq D \| \pi(\xi^2 \eta) \|_s^2, \quad 
\text{supp}(\gamma) \subset K_{j,m-1} \setminus K_{j,m-2}.
\label{eq:lemma3.5}
\end{equation}
Take $\theta = \xi^2 \eta - \gamma$. Again, $\pi(\theta) = 0$ and hence $\theta \in \widetilde{V}$. 
Therefore, we have 
\begin{equation}
a_Q(\psi_k^{(j)}, \theta) = 0.
\end{equation}
Additionally, $\text{supp}(\theta) \subset \Omega \setminus K_{j,m-2}$. 
Recall that, in \eqref{eq:lemma3.0}, we have $\text{supp}(\widetilde{\phi}_k^{(j)}) \subset K_j$.
Hence $\theta$ and $\widetilde{\phi}_k^{(j)}$ have disjoint supports, and
\begin{equation}
a_Q(\widetilde{\phi}_k^{(j)}, \theta) = 0.
\end{equation}
Therefore, we obtain
\begin{equation}
a_Q(\eta, \theta) = a_Q(\psi_k^{(j)}, \theta) - a_Q(\widetilde{\phi}_k^{(j)}, \theta) = 0.
\end{equation}
Note that $\xi^2 \eta = \theta + \gamma$. Using \eqref{eq:lemma3.5}, we have
\begin{equation}
\begin{split}
a_Q(\eta, \xi^2 \eta)
& = a_Q(\eta,\gamma) \\
& \leq \| \eta \|_{a_Q(K_{j,m-1} \setminus K_{j,m-2})} \| \gamma \|_{a_Q(K_{j,m-1} \setminus K_{j,m-2})} \\
& \leq D^{\frac{1}{2}} \| \eta \|_{a_Q(K_{j,m-1} \setminus K_{j,m-2})}  \| \pi(\xi^2 \eta) \|_{s(K_{j,m-1} \setminus K_{j,m-2})}.
\end{split}
\end{equation}
For any coarse element $K_{j'} \subset K_{j,m-1} \setminus K_{j,m-2}$, since $\pi(\eta) = 0$, we have
\begin{equation}
\| \pi(\xi^2 \eta) \|_{s(K_{j'})}
\leq \| \xi^2 \eta \|_{s(K_{j'})}
\leq \| \eta \|_{s(K_{j'})}
\leq \Lambda^{-\frac{1}{2}} \| \eta \|_{a_Q(K_{j'})}.
\end{equation}
Summing up over all $K_{j'} \subset K_{j,m-1} \setminus K_{j,m-2}$, we obtain
\begin{equation}
\| \pi(\xi^2 \eta) \|_{s(K_{j,m-1} \setminus K_{j,m-2})}
\leq \Lambda^{-\frac{1}{2}} \| \eta \|_{a_Q(K_{j,m-1} \setminus K_{j,m-2})}.
\end{equation}
Hence, the first term on the right hand side of \eqref{eq:lemma3.4} can be estimated by
\begin{equation}
a(\eta, \xi^2 \eta)
\leq D^{\frac{1}{2}} \Lambda^{-\frac{1}{2}} \| \eta \|^2_{a_Q(K_{j,m-1} \setminus K_{j,m-2})}.
\label{eq:lemma3.6a}
\end{equation}
For the second term on the right hand side of \eqref{eq:lemma3.4}, a similar argument gives
$\text{supp}(\xi \eta) \subset K_{j,m-1} \setminus K_{j,m-2}$, and
\begin{equation}
\| \xi \eta \|_s \leq \| \eta \|_{s(K_{j,m-1} \setminus K_{j,m-2})} \leq \Lambda^{-\frac{1}{2}} \| \eta \|_{a_Q(K_{j,m-1} \setminus K_{j,m-2})}.
\label{eq:lemma3.6b}
\end{equation}
Putting \eqref{eq:lemma3.4}, \eqref{eq:lemma3.6a} and \eqref{eq:lemma3.6b} together, we have
\begin{equation}
\| \eta \|^2_{a_Q(\Omega \setminus K_{j,m-1})} 
\leq (2 + D^\frac{1}{2}) \Lambda^{-\frac{1}{2}} \| \eta \|^2_{a_Q(K_{j,m-1} \setminus K_{j,m-2})}
\leq 3D^\frac{1}{2} \Lambda^{-\frac{1}{2}} \| \eta \|^2_{a_Q(K_{j,m-1} \setminus K_{j,m-2})}.
\end{equation}
Therefore,
\begin{equation}
\| \eta \|^2_{a_Q(\Omega \setminus K_{j,m-2})} 
= \| \eta \|^2_{a_Q(\Omega \setminus K_{j,m-1})} + \| \eta \|^2_{a_Q( K_{j,m-1} \setminus K_{j,m-2})}
\geq \left(1 + \dfrac{\Lambda^\frac{1}{2}}{3D^\frac{1}{2}}\right) \| \eta \|^2_{a_Q(\Omega \setminus K_{j,m-1})}.
\end{equation}
Inductively, we have
\begin{equation}
\| \eta \|^2_{a_Q(\Omega \setminus K_{j,m-1})} 
\leq \left(1 + \dfrac{\Lambda^\frac{1}{2}}{3D^\frac{1}{2}}\right)^{1-m} \| \eta \|^2_{a_Q(\Omega \setminus K_{j})}
\leq \left(1 + \dfrac{\Lambda^\frac{1}{2}}{3D^\frac{1}{2}}\right)^{1-m} \| \eta \|^2_{a_Q}.
\label{eq:lemma3.7}
\end{equation}
Finally, by the energy minimzing property of $\psi_k^{(j)}$ and \eqref{eq:lemma3.0}, we have
\begin{equation}
\| \eta \|_{a_Q} = \| \psi_k^{(j)} - \widetilde{\phi}_k^{(j)} \|_{a_Q} \leq 2 \| \widetilde{\phi}_k^{(j)} \|_{a_Q} \leq 2 D^\frac{1}{2}  \| \phi_k^{(j)} \|_{s(K_j)}.
\label{eq:lemma3.8}
\end{equation}
Combining \eqref{eq:lemma3.3}, \eqref{eq:lemma3.7} and \eqref{eq:lemma3.8}, we obtain our desired result.
\end{proof}
\end{lemma}

%We will need one more estimate in order to prove the usefulness of the localized multiscale basis functions.
%\begin{lemma}
%\label{lemma4}
%With the same notations and assumptions in Lemma~\ref{lemma3}, we have
%\begin{equation}
%\left\| \sum_{j=1}^N (\psi_k^{(j)} - \psi_{k,ms}^{(j)}) \right\|_{a_Q}^2 \leq C(m+1)^d \sum_{j=1}^N  \left\| \psi_k^{(j)} - \psi_{k,ms}^{(j)} \right\|_{a_Q}^2.
%\end{equation}
%\begin{proof}
%\end{proof}
%\end{lemma}

The above lemma motivates us to define localized multiscale basis functions in \eqref{eq:min1}. 
The following lemma suggests that, 
similar to the projection operator $R_{glo}$ onto the global multiscale finite element space, 
the projection operator $R_{ms}$ onto our localized multiscale finite element space 
also has a good approximation property with respect to the $a_Q$-norm and $L^2$-norm.
\begin{lemma}
\label{lemma4}
Let $u \in D(\mathcal{A})$. 
Let $m \geq 2$ be the number of coarse grid layers in the oversampling regions in \eqref{eq:min1}. 
If $m = O\left(\log\left(\dfrac{\overline{\kappa}}{H}\right)\right)$, then we have 
\begin{equation}
\| u - R_{ms} u \|_{a_Q} \leq CH \underline{\kappa}^{-\frac{1}{2}} \Lambda^{-\frac{1}{2}} \| \mathcal{A} u \|_{[L^2(\Omega)]^2},
\label{eq:a_approx_ms}
\end{equation}
and
\begin{equation}
\| u - R_{ms} u \|_{[L^2(\Omega)]^2} \leq CH^2 \underline{\kappa}^{-1} \Lambda^{-1} \| \mathcal{A} u \|_{[L^2(\Omega)]^2}.
%\| u - R_{ms} u \|_c \leq CH^2 \Lambda^{-1} \| \mathcal{A} u \|_{[L^2(\Omega)]^2}.
\label{eq:c_approx_ms}
\end{equation}
\begin{proof}
We write $R_{glo} u = \sum_{j=1}^N \sum_{k=1}^{L_j} \alpha_k^{(j)} \psi_k^{(j)}$, and 
define $w = \sum_{j=1}^N \sum_{k=1}^{L_j} \alpha_k^{(j)} \psi_{k,ms}^{(j)} \in V_{ms}$. 
By the Galerkin orthogonality in \eqref{eq:elliptic_proj_ms}, we have
\begin{equation}
\| u - R_{ms} u \|_{a_Q} \leq \| u - w \|_{a_Q} \leq \| u - R_{glo} u\|_{a_Q} + \| R_{glo} u - w \|_{a_Q}.
\label{eq:lemma5.1}
\end{equation}
%Using Lemma~\ref{lemma3} and Lemma~\ref{lemma4}, we see that
Using Lemma~\ref{lemma3}, we see that
\begin{equation}
\begin{split}
\| R_{glo} u - w \|_{a_Q}^2
& = \left\| \sum_{j=1}^N \sum_{k=1}^{L_j} \alpha_k^{(j)} (\psi_k^{(j)} - \psi_{k,ms}^{(j)}) \right\|_{a_Q}^2 \\
& \leq C(m+1)^d \sum_{j=1}^N \left\| \sum_{k=1}^{L_j} \alpha_k^{(j)} (\psi_k^{(j)} - \psi_{k,ms}^{(j)}) \right\|_{a_Q}^2 \\
& \leq CE(m+1)^d \sum_{j=1}^N \left\| \sum_{k=1}^{L_j} \alpha_k^{(j)} \phi_k^{(j)} \right\|_s^2 \\
& = CE(m+1)^d \| R_{glo} u \|_s^2,
\end{split}
\label{eq:lemma5.2}
\end{equation}
where the last equality follows from the orthogonality of the eigenfunctions in \eqref{eq:spectral_prob}. 
Combining \eqref{eq:lemma5.1}, \eqref{eq:lemma5.2}, together with \eqref{eq:a_approx_glo} in Lemma~\ref{thm:elliptic-glo}, we have
\begin{equation}
\| u - R_{ms} u \|_{a_Q} \leq CH \underline{\kappa}^{-\frac{1}{2}} \Lambda^{-\frac{1}{2}} \| \mathcal{A} u \|_{[L^2(\Omega)]^2} + CE^\frac{1}{2}(m+1)^{\frac{d}{2}} \| R_{glo} u \|_s.
\end{equation}
Next, we are going to estimate $\| R_{glo} u \|_s$. Using the fact that $\vert \nabla \chi_k \vert = O(H^{-1})$, we have
\begin{equation}
\| R_{glo} u \|_s^2 \leq CH^{-2} \overline{\kappa} \| R_{glo} u \|_{[L^2(\Omega)]^2}^2.
\end{equation}
Then, by Poincar\'{e} inequality, we have
\begin{equation}
\| R_{glo} u \|^2_{[L^2(\Omega)]^2} \leq C_p \underline{\kappa}^{-1} \| R_{glo} u \|_{a_Q}^2.
\end{equation}
By taking $v = R_{glo} u$ in \eqref{eq:elliptic_proj}, we obtain
\begin{equation}
\| R_{glo} u \|_{a_Q}^2 = a_Q(u, R_{glo}u) = (\mathcal{A}u, R_{glo}u) \leq CH\underline{\kappa}^{-\frac{1}{2}} \| \mathcal{A}u\|_{[L^2(\Omega)]^2} \| R_{glo}u \|_s.
\end{equation}
Combining these estimates, we have
\begin{equation}
\| R_{glo} u \|_s \leq CH^{-1} \overline{\kappa} \underline{\kappa}^{-\frac{1}{2}} \| \mathcal{A}u\|_{[L^2(\Omega)]^2}.
\end{equation}
To obtain our desired result, we need
\begin{equation}
H^{-2}\overline{\kappa}(m+1)^{\frac{d}{2}}E^\frac{1}{2} = O(1).
\end{equation}
Taking logarithm, we have
\begin{equation}
\log(H^{-2}) + \log(\overline{\kappa}) + \dfrac{d}{2} \log(m+1) + \dfrac{1-m}{2} \log\left(1+\dfrac{\Lambda^\frac{1}{2}}{3D^\frac{1}{2}}\right) = O(1).
\end{equation}
Thus, taking $m = O\left(\log\left(\dfrac{\overline{\kappa}}{H}\right)\right)$ completes the proof of \eqref{eq:a_approx_ms}. 
The proof of \eqref{eq:c_approx_ms} follows from a duality argument as in Lemma~\ref{thm:elliptic-glo}.
\end{proof}
\label{thm:elliptic-ms}
\end{lemma}

We are now ready to establish our main theorem, which estimates the error between the solution $p$ and the multiscale solution $p_{ms}$.
\begin{theorem}
\label{theorem1}
Suppose $f \in [L^2(\Omega)]^2$. 
Let $m \geq 2$ be the number of coarse grid layers in the oversampling regions in \eqref{eq:min1}.
Let $p$ be the solution of \eqref{eq:sol_weak} and
$p_{ms}$ be the solution of \eqref{eq:sol_ms}.
If $m = O\left(\log\left(\dfrac{\overline{\kappa}}{H}\right)\right)$, then we have 
\begin{equation}
\| p(T,\cdot) - p_{ms}(T,\cdot) \|_c^2 + \int_0^T \| p - p_{ms} \|_{a_Q}^2 dt \leq CH^2 \underline{\kappa}^{-1} \Lambda^{-1} \left(\| p^0 \|_{a_Q}^2 + \int_0^T \| f \|_{[L^2(\Omega)]^2}^2 \, dt\right).
\end{equation}
\begin{proof}
Taking $v = \dfrac{\partial p}{\partial t}$ in \eqref{eq:sol_weak}, we have
\begin{equation}
\left\| \dfrac{\partial p}{\partial t} \right\|_c^2 + \dfrac{1}{2} \dfrac{d}{dt} \| p \|_{a_Q}^2 = \left( f, \dfrac{\partial p}{\partial t} \right) \leq C \| f \|_{[L^2(\Omega)]^2}^2 + \dfrac{1}{2} \left\| \dfrac{\partial p}{\partial t} \right\|_c^2.
\end{equation}
Integrating over $(0,T)$, we have
\begin{equation}
\dfrac{1}{2} \int_0^T \left\| \dfrac{\partial p}{\partial t} \right\|_c^2 dt + \dfrac{1}{2} \| p(T, \cdot) \|_{a_Q}^2 
\leq C \left(\| p^0 \|_{a_Q}^2 + \int_0^T \| f \|_{[L^2(\Omega)]^2}^2 dt\right).
\label{eq:thm1.1}
\end{equation}
Similarly, taking $v = \dfrac{\partial p_{ms}}{\partial t}$ in \eqref{eq:sol_ms} and integrating over $(0,T)$, we have
\begin{equation}
\dfrac{1}{2} \int_0^T \left\| \dfrac{\partial p_{ms}}{\partial t} \right\|_c^2 dt + \dfrac{1}{2} \| p_{ms}(T, \cdot) \|_{a_Q}^2 
\leq C \left( \| p^0 \|_{a_Q}^2 + \int_0^T \| f \|_{[L^2(\Omega)]^2}^2 dt\right).
\label{eq:thm1.2}
\end{equation}
On the other hand, from $\eqref{eq:sol_weak}$, we see that
\begin{equation}
\mathcal{A}p = f - \mathcal{C} \dfrac{\partial p}{\partial t},
\end{equation}
and therefore
\begin{equation}
\| \mathcal{A}p \|_{[L^2(\Omega)]^2} \leq C \left(\| f \|_{[L^2(\Omega)]^2} + \left\| \dfrac{\partial p}{\partial t} \right\|_c\right).
\label{eq:thm1.3}
\end{equation}
By the definition of $p$ in \eqref{eq:sol_weak} and $p_{ms}$ in \eqref{eq:sol_ms}, for all $v \in V_{ms}, t \in (0,T)$, we have
\begin{equation}
c \left(\dfrac{\partial (p - p_{ms})}{\partial t}, v \right) + a_Q(p - p_{ms}, v) = 0.
\end{equation}
Therefore, we have
\begin{equation}
\begin{split}
& \dfrac{1}{2} \dfrac{d}{dt} \| p - p_{ms} \|_c^2 + \| p - p_{ms} \|_{a_Q}^2 \\
& = c \left(\dfrac{\partial (p - p_{ms})}{\partial t}, p - p_{ms} \right) + a_Q(p - p_{ms}, p - p_{ms})  \\
& = c \left(\dfrac{\partial (p - p_{ms})}{\partial t}, p - R_{ms}p \right) + a_Q(p - p_{ms}, p - R_{ms}p )  \\
& \leq \left\| \dfrac{\partial (p - p_{ms})}{\partial t} \right\|_c \| p - R_{ms}p \|_c + \|p - p_{ms}\|_{a_Q} \| p - R_{ms} p \|_{a_Q} \\
& \leq \left(\left\| \dfrac{\partial p}{\partial t} \right\|_c + \left\| \dfrac{\partial p_{ms}}{\partial t} \right\|_c\right) \| p - R_{ms}p \|_c + \dfrac{1}{2} \|p - p_{ms}\|_{a_Q}^2 + \dfrac{1}{2} \| p - R_{ms} p \|_{a_Q}^2.
\end{split}
\end{equation}
Integrating over $(0,T)$ and using \eqref{eq:thm1.3} with Lemma~\ref{thm:elliptic-ms}, we have
\begin{equation}
\begin{split}
& \dfrac{1}{2} \| p(T,\cdot) - p_{ms}(T,\cdot) \|_c^2 + \dfrac{1}{2} \int_0^T \| p - p_{ms} \|_{a_Q}^2 dt \\
& \leq \int_0^T \left(\left\| \dfrac{\partial p}{\partial t} \right\|_c + \left\| \dfrac{\partial p_{ms}}{\partial t} \right\|_c\right) \| p - R_{ms}p \|_c dt + \dfrac{1}{2}\int_0^T \| p - R_{ms} p \|_{a_Q}^2 dt \\
& \leq \left( \int_0^T \left(\left\| \dfrac{\partial p}{\partial t} \right\|_c + \left\| \dfrac{\partial p_{ms}}{\partial t} \right\|_c\right)^2 dt \right)^\frac{1}{2} \left(\int_0^T \| p - R_{ms}p \|_c^2 \, dt \right)^\frac{1}{2} + \dfrac{1}{2}\int_0^T \| p - R_{ms} p \|_{a_Q}^2 dt \\
& \leq \left( \int_0^T \left(\left\| \dfrac{\partial p}{\partial t} \right\|_c + \left\| \dfrac{\partial p_{ms}}{\partial t} \right\|_c\right)^2 dt \right)^\frac{1}{2} \left(\int_0^T C H^4 \underline{\kappa}^{-2} \Lambda^{-2}  \left(\| f \|_{[L^2(\Omega)]^2} + \left \| \dfrac{\partial p}{\partial t} \right\|_c \right)^2 \, dt \right)^\frac{1}{2} + \\
& \quad \quad \int_0^T C H^2 \underline{\kappa}^{-1} \Lambda^{-1}  \left(\| f \|_{[L^2(\Omega)]^2} + \left \| \dfrac{\partial p}{\partial t} \right\|_c \right)^2 \, dt \\
& \leq CH^2 \underline{\kappa}^{-1} \Lambda^{-1} \int_0^T \left( \left\| \dfrac{\partial p}{\partial t} \right\|_c^2 +  \left\| \dfrac{\partial p_{ms}}{\partial t} \right\|_c^2 + \| f \|_{[L^2(\Omega)]^2}^2 \right) dt.
\end{split}
\label{eq:thm1.4}
\end{equation}
Finally, combining \eqref{eq:thm1.1}, \eqref{eq:thm1.2} and \eqref{eq:thm1.4}, we obtain our desired result.
\end{proof}
\end{theorem}

\section{Numerical Examples}\label{sec:numerical}

\begin{figure}[ht!]
\centering
\includegraphics[width=0.48\linewidth]{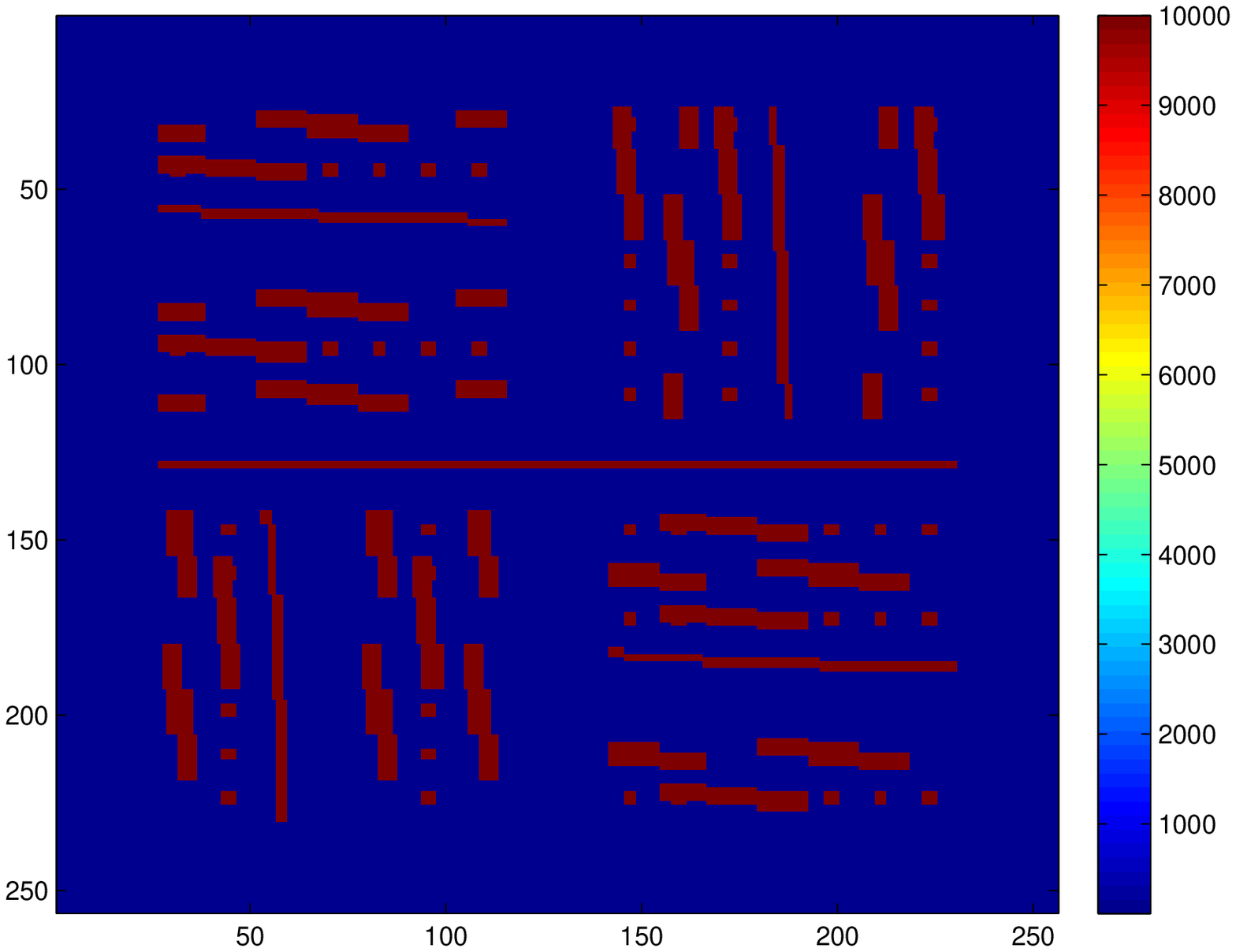}
\includegraphics[width=0.48\linewidth]{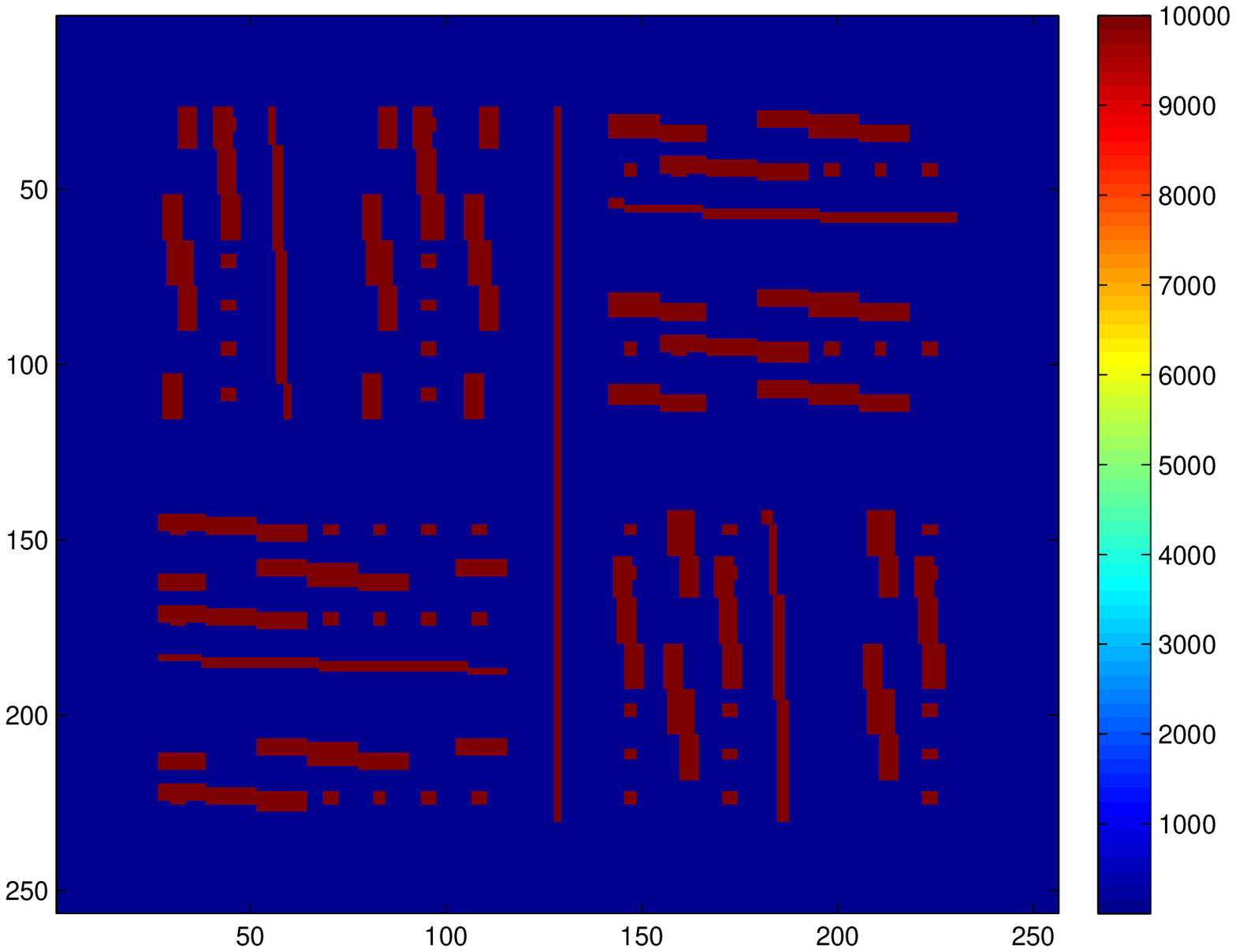}
\caption{Media used in numerical experiments. $\kappa_1$ (left) and $\kappa_2$ (right).}
\label{fig:kappa1}
\end{figure}

In this section, we present two numerical examples. 
We perform numerical experiments with high-contrast media 
to see the orders of convergence 
of our proposed method in energy norm and $L^2$ norm. 
We will also study the effects of the number of oversampling layers 
$m$ on the quality of the approximations.  
In all the experiments, we take the spatial domain to be $\Omega = (0,1)^2$ and 
the fine mesh size to be $h = 1/256$. 
An example of the media $\kappa_1$ and $\kappa_2$ used in the experiments 
is illustrated in FIgure~\ref{fig:kappa1}. 
In the figure, the contrast values, i.e. the ratio of the maximum and the minimum in $\Omega$, of the media are
$\overline{\kappa}_1 = 10^4$ and $\overline{\kappa}_2 = 10^4$. 
We will also see the effects of the contrast values of the media on the error, 
while the configurations of the media remain unchanged. 

\subsection{Experiment 1}

In this experiment, we consider the dual continuum model in the steady state, i.e. 
\begin{equation}
\begin{split}
- \text{div}(\kappa_1 \nabla p_1) + \rho \sigma(p_1 - p_2) = \rho f_1, \\
- \text{div}(\kappa_2 \nabla p_2) - \rho \sigma(p_1 - p_2) = \rho f_2,
\end{split}
\label{eq:steady}
\end{equation}
where the configuration of the media $\kappa_1$ and $\kappa_2$ 
are illustrated in FIgure~\ref{fig:kappa1}. 
The conductivity values in the background are fixed to be $\kappa_{1,m} = 1$ and $\kappa_{2,m} = 1$, 
while the conductivity values $\kappa_{1,f}$ and $\kappa_{2,f}$ in the channels are high.
The physical constants are set to be $\rho = 1$ and $\sigma = 1$. 
The source functions are taken as $f_1(x,y) =  2\pi^2\sin(\pi x) \sin(\pi y)$ 
and $f_2(x,y) = 1$ for all $(x,y) \in \Omega$. 
The steady-state equation \eqref{eq:steady} has a weak formulation: 
find $p = (p_1, p_2)$ with $p_i \in V$ such that
\begin{equation}
a_Q(p,v) = (f,v),
\end{equation}
for all $v = (v_1, v_2)$ with $v_i \in V$. 
The numerical solution is then given by: find $p_{ms} = (p_{ms,1}, p_{ms,2})$ 
with $p_{ms,i} \in V_{ms}$ such that  
\begin{equation}
a_Q(p_{ms},v) = (f,v),
\end{equation}
for all $v = (v_1, v_2)$ with $v_i \in V_{ms}$. 
In other words, we have $p_{ms} = R_{ms} p$ according to the definition \eqref{eq:elliptic_proj_ms}, 
and the theoretical orders of convergence follow Lemma~\ref{lemma4}. 

\begin{figure}[ht!]
\centering
\includegraphics[width=0.48\linewidth]{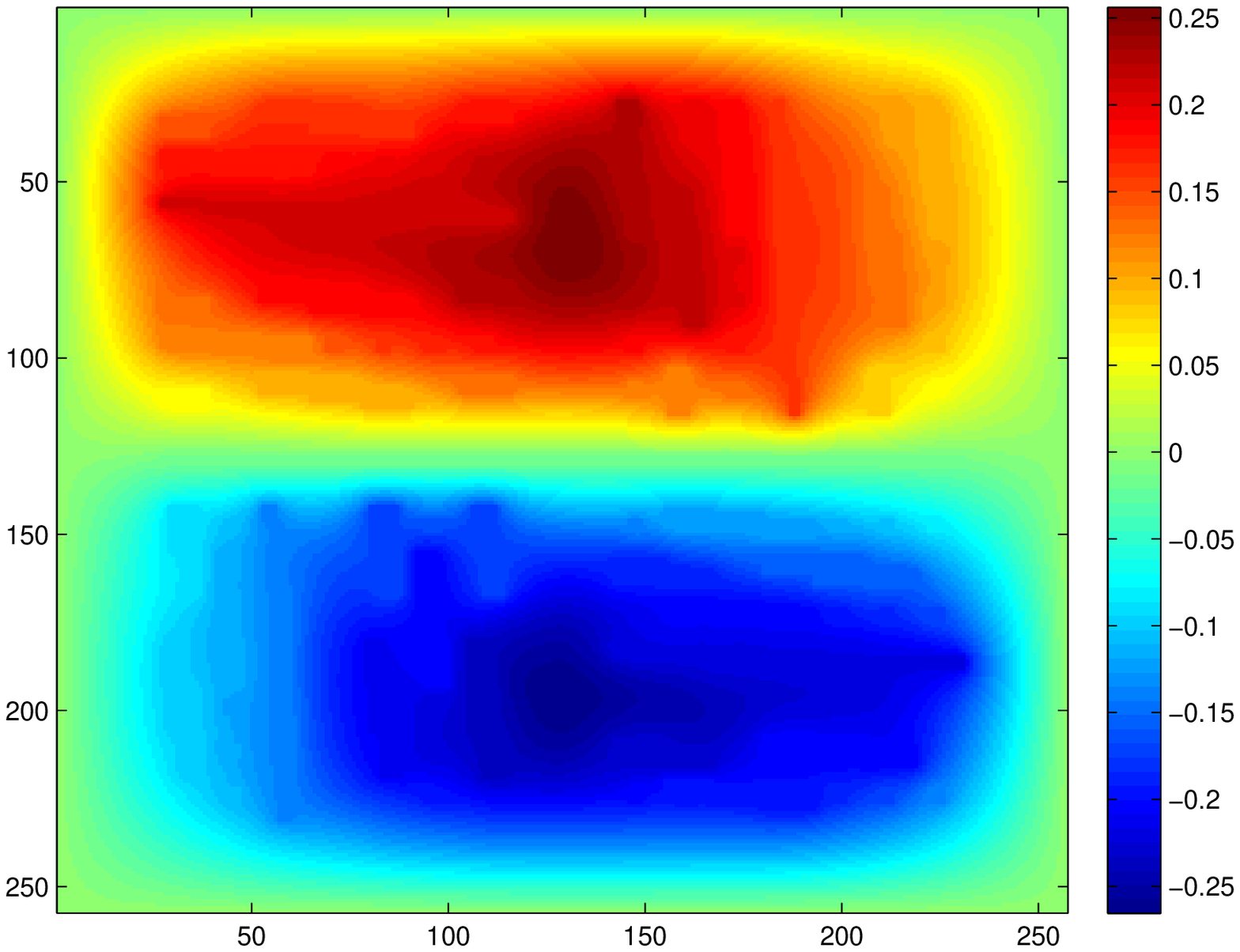}
\includegraphics[width=0.48\linewidth]{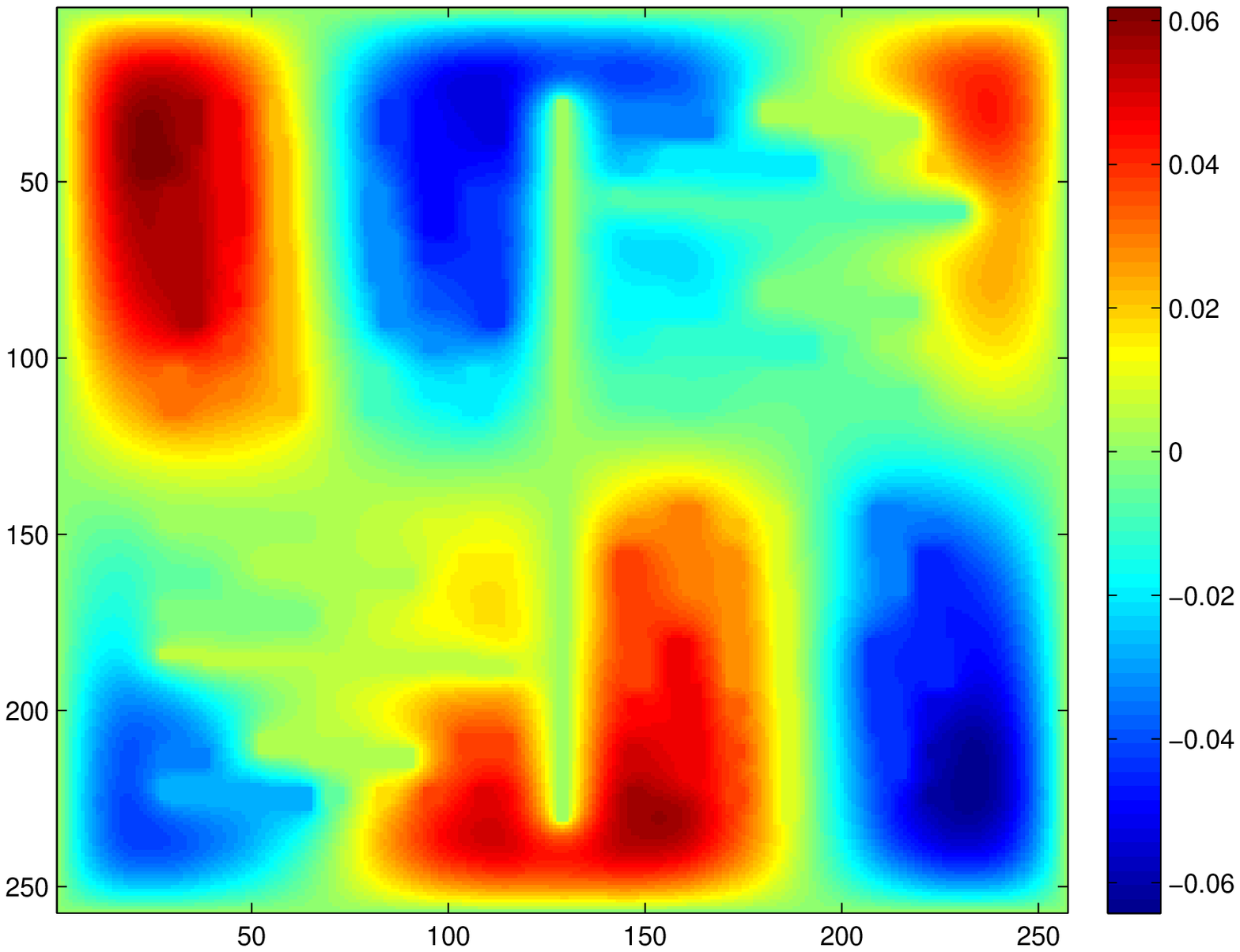}
\caption{Plots of numerical solution: $p_{ms,1}$ (left) and $p_{ms,2}$ (right)}
\label{fig:exp1}
\end{figure}

Figure~\ref{fig:exp1} illustrates the numerical solution of the steady-state flow problem. 
Tables~\ref{tab:exp1.1}--\ref{tab:exp1.3} record 
the error in $L^2$ norm and $a_Q$ norm with various settings. 
In Table~\ref{tab:exp1.1}, we take the conductivity values in the channels to be 
$\kappa_{1,f} = 10^4$ and $\kappa_{2,f} = 10^6$. 
We use $6$ basis functions per oversampled region since 
there are $6$ small eigenvalues in the spectrum, 
and according to our analysis,  
we need to include the first $6$ spectral basis functions 
in the auxiliary space to have good convergence. 
As we refine coarse mesh size $H$, 
we fix the number of oversampling layers to be $m \approx 9 \log(1/H) / \log(64)$, 
which is suggested by our analysis. 
The results show that the numerical approximations are very accurate, 
and the errors converge with refinement of the coarse mesh size. 
Table~\ref{tab:exp1.2} shows the same quantities when we reduce 
the number of basis functions used in each coarse region is reduced to $4$. 
By comparing to Table~\ref{tab:exp1.1}, it can be seen that the errors are larger than 
those when we use $6$ basis functions. 
Table~\ref{tab:exp1.3} compares the $a_Q$ error with various combinations of 
number of layers $m$ and contrast value $\overline{\kappa}$, 
where the conductivity values in the channels are the same, 
with $6$ basis functions per coarse region and coarse mesh size $H = 1/16$. 
It can be seen that with a larger oversampled region, the error increases. 
On the other hand, the error increases with the contrast value.

\begin{table}[ht!]
\centering
\begin{tabular}{|c|c||c|c||c|c|}
\hline
$H$ & $m$ & $a_Q$ error & order & $L^2$ error & order \\
\hline
$1/8$ & 4 & 33.4293\% & -- & 15.8783\% & -- \\
$1/16$ & 6 & 5.7191\% & 2.55 & 0.6265\% & 4.66 \\
$1/32$ & 7 & 1.2437\% & 2.20 & 0.0504\% & 3.64 \\
$1/64$ & 9 & 0.3585\% & 1.79 & 0.0067\% & 2.91 \\
\hline
\end{tabular}
\caption{History of convergence with $6$ basis functions in Experiment 1.}
\label{tab:exp1.1}
\end{table}

\begin{table}[ht!]
\centering
\begin{tabular}{|c|c||c|c||c|c|}
\hline
$H$ & $m$ & $a_Q$ error & order & $L^2$ error & order \\
\hline
$1/8$ & 4 & 43.9247\% & -- & 34.2923\% & -- \\
$1/16$ & 6 & 7.7963\% & 2.49 & 1.0463\% & 5.03 \\
$1/32$ & 7 & 1.5417\% & 2.34 & 0.0709\% & 3.88 \\
$1/64$ & 9 & 0.4993\% & 1.63 & 0.0124\% & 2.52 \\
\hline
\end{tabular}
\caption{History of convergence with $4$ basis functions in Experiment 1.}
\label{tab:exp1.2}
\end{table}

\begin{table}[ht!]
\centering
\begin{tabular}{|c||c|c|c|}
\hline
$m$ & $\overline{\kappa} = 10^4$ & $\overline{\kappa} = 10^5$ & $\overline{\kappa} = 10^6$ \\
\hline
3 & 22.4683\% & 51.0835\% & 69.4279\% \\
4 & 6.3274\% & 10.1892\% & 25.6786\% \\
5 & 5.7205\% & 5.7978\% & 6.4329\% \\
6 & 5.7122\% & 5.7220\% & 5.7231\% \\
\hline
\end{tabular}
\caption{Comparison of $a_Q$ error with different number of layers $m$ and contrast value $\overline{\kappa}$ in Experiment 1.}
\label{tab:exp1.3}
\end{table}

\subsection{Experiment 2}

\begin{figure}[ht!]
\centering
\includegraphics[width=0.6\linewidth]{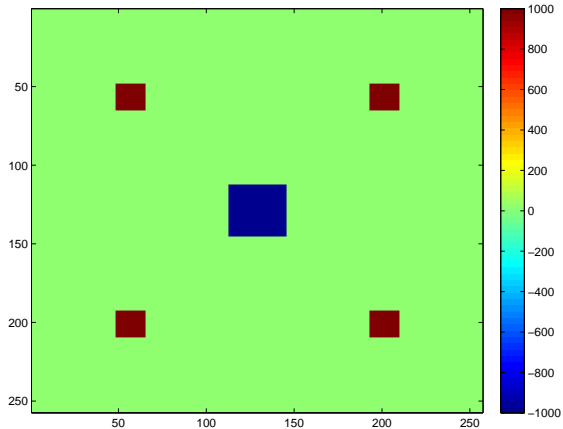}
\caption{Source function $f_2$ in Experiment 2.}
\label{fig:source2}
\end{figure}

In this experiment, we consider the time-dependent dual continuum model \eqref{eq:dc}. 
We are interested in finding a numerical approximation in
the temporal domain $[0,T]$, where the final time is set to be $T = 5$. 
The configuration of the media $\kappa_1$ and $\kappa_2$ 
are illustrated in FIgure~\ref{fig:kappa1}. 
The conductivity values in the background are set to be 
$\kappa_{1,m} = 10^{-1}$ and $\kappa_{2,m} = 10^0$, 
while the values in the channels are taken as 
$\kappa_{1,f} = 10^4$ and $\kappa_{2,f} = 10^6$. 
The velocities in the background are taken as 
$c_{1,m} = 10^{1}$ and $c_{2,m} = 10^3$, 
while the values in the channels are taken as 
$c_{1,f} = 10^2$ and $c_{2,f} = 10^4$. 
The physical constants are set to be $\rho = 1$ and $\sigma = 25$. 
The source functions are taken as time-independent, 
where $f_1(t,x,y) = 0$ for all $(t,x,y) \in [0,T] \times \Omega$
and $f_2$ is depicted in Figure~\ref{fig:source2}. 
The initial condition is given as $p_1(0,x,y) = 0$ and $p_2(0,x,y) = 0$ for all $(x,y) \in \Omega$. 

\begin{figure}[ht!]
\centering
\includegraphics[width=0.48\linewidth]{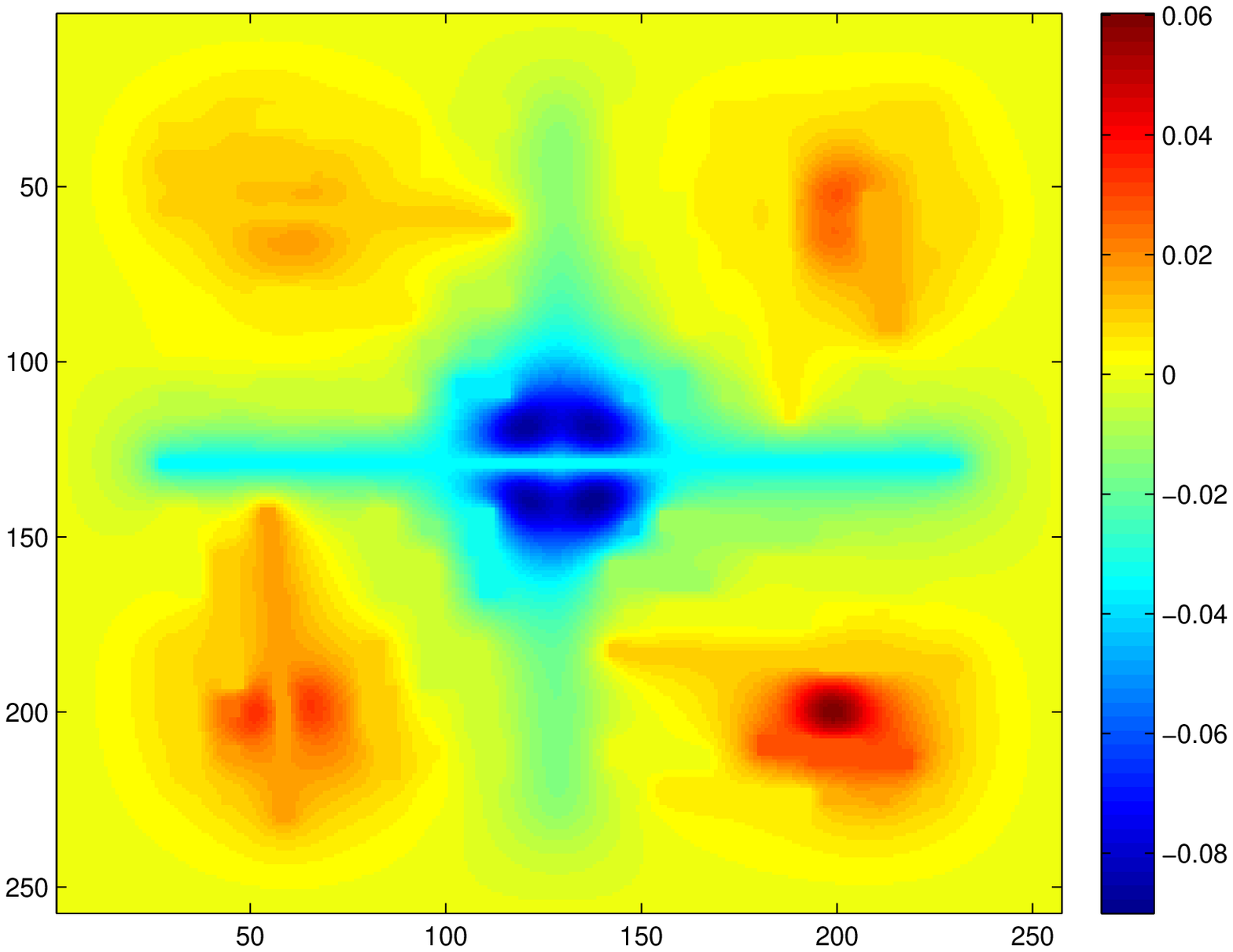}
\includegraphics[width=0.48\linewidth]{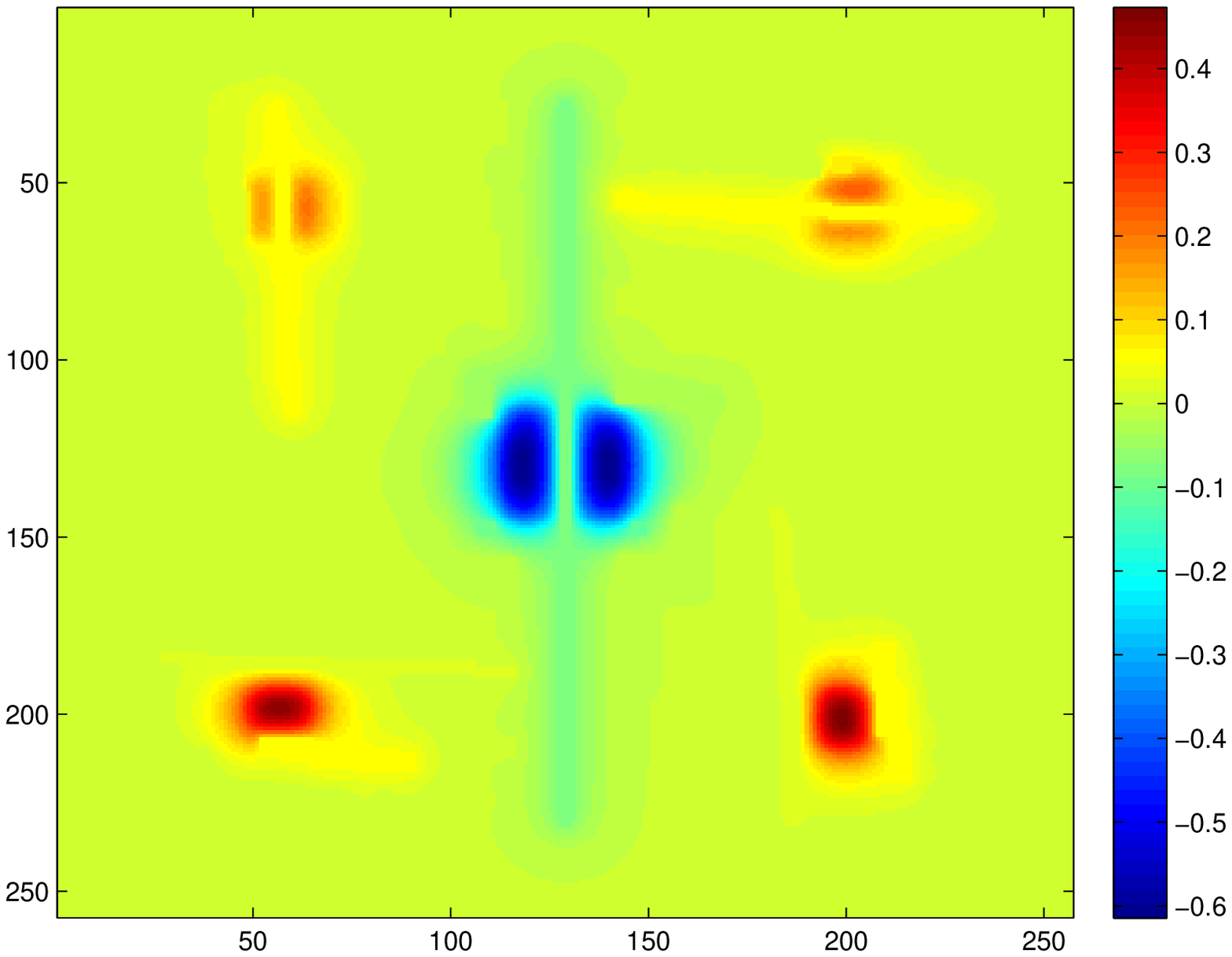}\\
\includegraphics[width=0.48\linewidth]{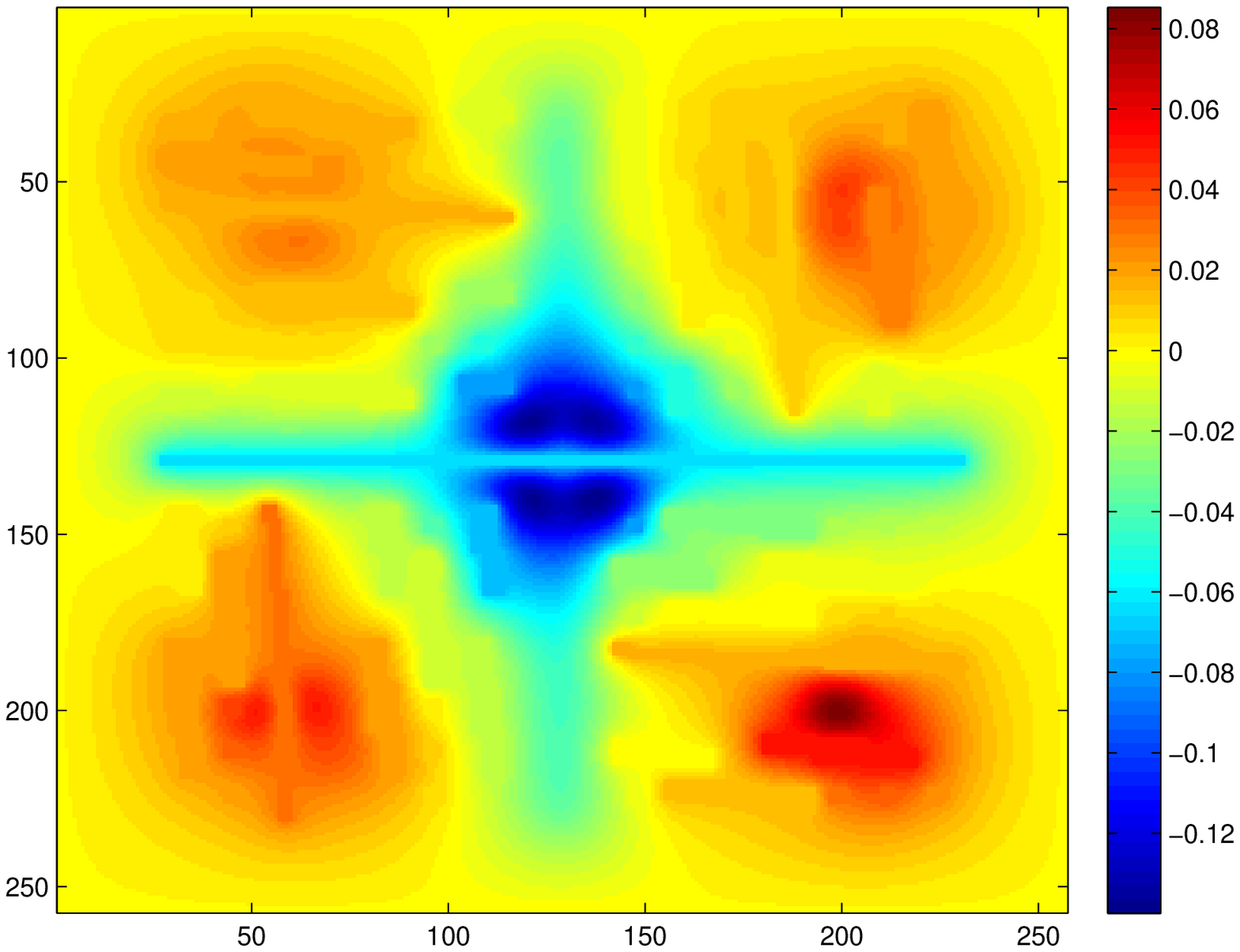}
\includegraphics[width=0.48\linewidth]{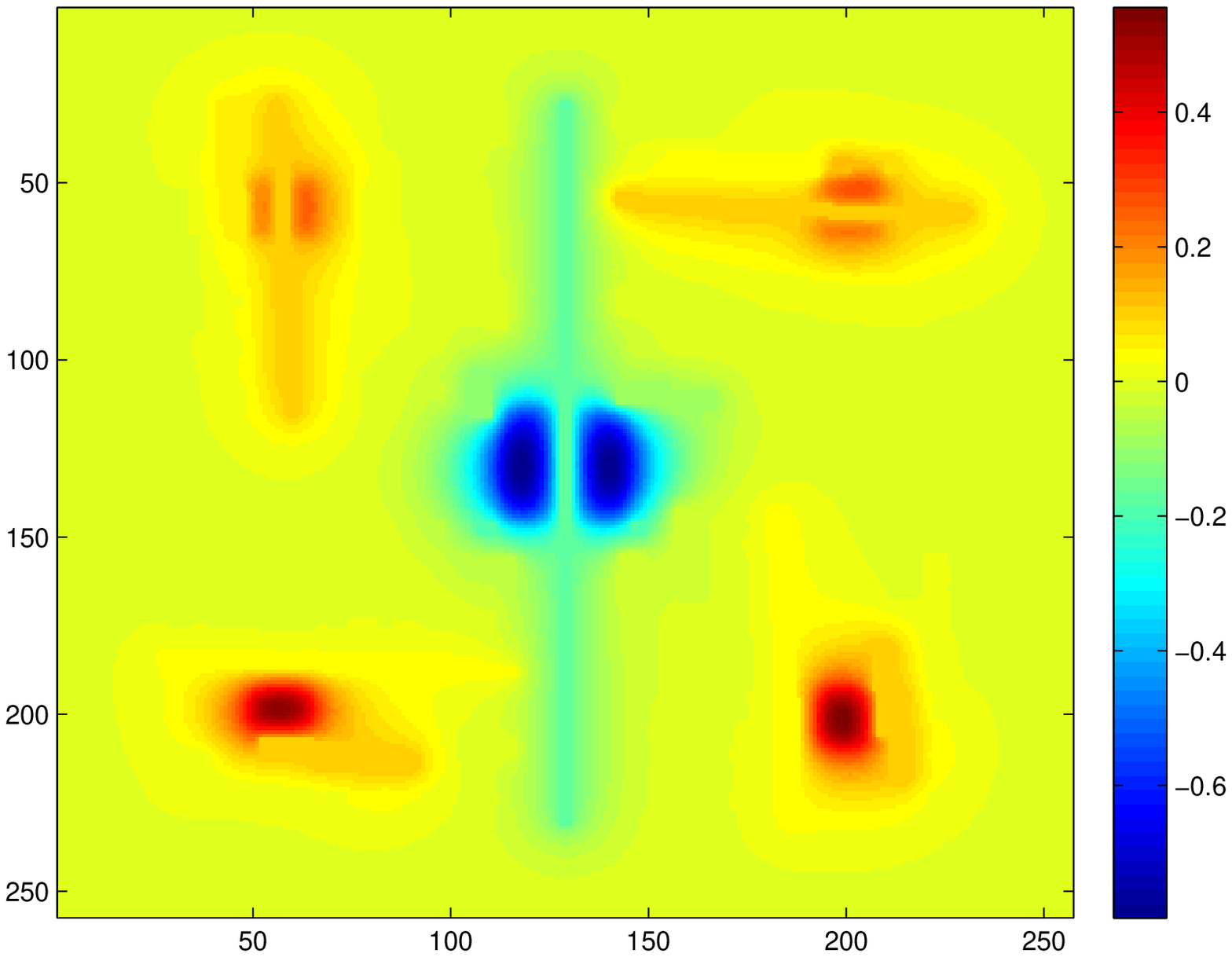}\\
\includegraphics[width=0.48\linewidth]{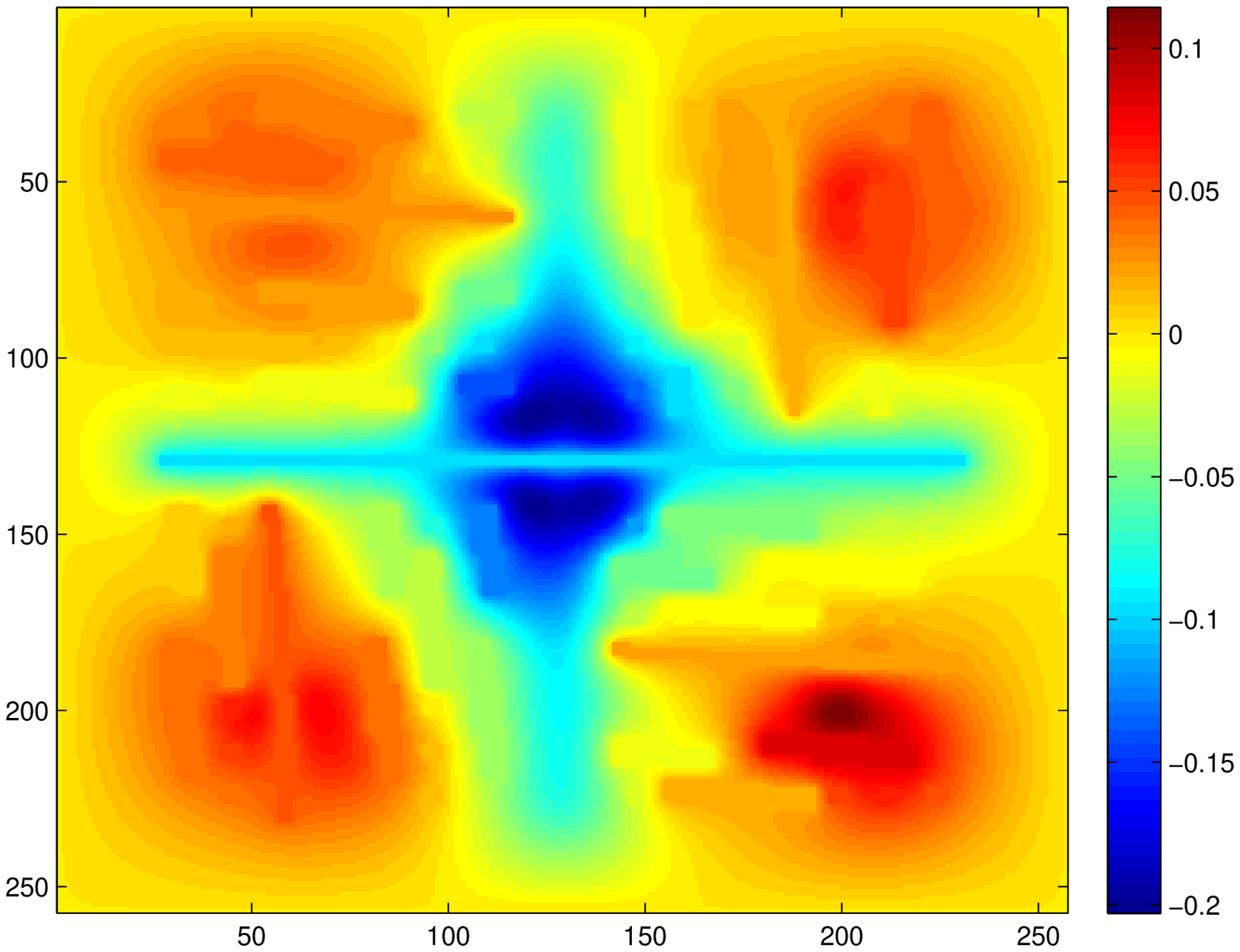}
\includegraphics[width=0.48\linewidth]{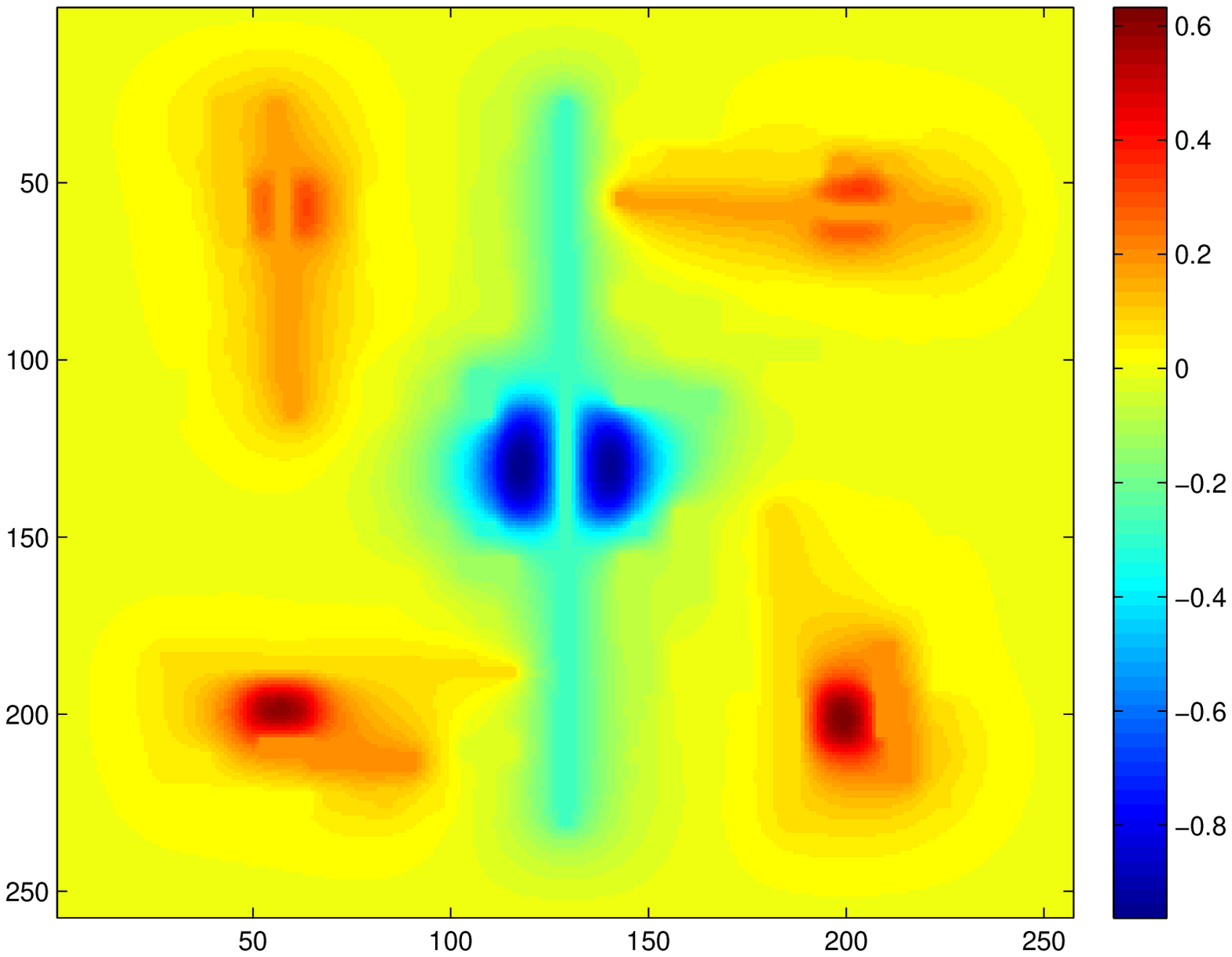}
\caption{Plots of numerical solution at different time instants: $p_{ms,1}$ (left) and $p_{ms,2}$ (right) in Experiment 2.}
\label{fig:exp2}
\end{figure}

Figure~\ref{fig:exp2} illustrates the numerical solutions 
at time instants $t = 1.25$, $t = 2.5$ and $t = 5$ respectively. 
Tables~\ref{tab:exp2} records
the error in $L^2$ norm and $a_Q$ norm with $6$ basis functions per oversampled region
and number of oversampling layers set to be $m \approx 9 \log(1/H) / \log(64)$. 
Again, the results show that the numerical approximations are very accurate, 
and the errors converge with with refinement of the coarse mesh size. 

% Exp 11

\begin{table}[ht!]
\centering
\begin{tabular}{|c|c|c||c|c||c|c|}
\hline
$H$ & $m$ & $\Delta t$ & $a_Q$ error & order & $L^2$ error & order \\
\hline
$1/8$ & 4 & 1 & 92.0441\% & -- & 58.6453\% & -- \\
$1/16$ & 6 & 0.5 & 20.9725\% & 2.13 & 5.2984\% & 3.47 \\
$1/32$ & 7 & 0.25 & 6.7504\% & 1.64 & 0.7718\% & 2.78 \\
$1/64$ & 9 & 0.125 & 1.9074\% & 1.82 & 0.0934\% & 3.05 \\
\hline
\end{tabular}
\caption{History of convergence with $6$ basis functions in Experiment 2.}
\label{tab:exp2}
\end{table}

\section{Conclusions}\label{sec:conclusions}

In this paper, we present the CEM-GMsFEM for a dual continuum model.
The method is based on a set of multiscale basis functions.
To find the basis, we first obtain the auxiliary basis functions, which are important to identify high contrast channels and fracture networks. 
Then, we solve an energy minimization with some constraints related to the auxiliary functions.
We show that the basis functions are localized and that the resulting method has a mesh dependent convergence.
Finally, we present some numerical results to confirm the theory.

\section*{Acknowledgements}

EC's work is partially supported by Hong Kong RGC General Research Fund (Project 14304217)
and CUHK Direct Grant for Research 2017-18.

\bibliographystyle{plain}
\bibliography{references2,references1,references,references_outline}

\end{document}